\documentclass[a4paper, 12pt, 
oneside, onecolumn]{elsarticle}
\usepackage{amssymb, amsthm}
\usepackage{verbatim}
\usepackage{graphicx}
\usepackage{subcaption}
\usepackage{float}    
\usepackage{amsthm}
\usepackage{amsmath, physics,cleveref}
\usepackage{mathrsfs}
\usepackage[misc]{ifsym}
\usepackage{diagbox}
\usepackage{graphicx} 
\usepackage{amsmath,latexsym,amssymb,amsfonts,amsbsy}
\usepackage{algorithmic}
\usepackage{algorithm}
\usepackage{color,xcolor}
\usepackage{booktabs}
\usepackage{threeparttable}
\usepackage{multicol}
\usepackage{epstopdf}
\usepackage{multirow}
\numberwithin{equation}{section}
\newtheorem{theorem}{Theorem}[section]
\newtheorem{lemma}[theorem]{Lemma}

\newtheorem{remark}[theorem]{Remark}

\newtheorem{corollary}[theorem]{Corollary}

\begin{document} 
\makeatletter
\def\ps@pprintTitle{%
  \let\@oddhead\@empty
  \let\@evenhead\@empty
  \let\@oddfoot\@empty
  \let\@evenfoot\@empty
}
\makeatother
\begin{frontmatter}
    \title{Data assimilation with model errors}
    
\author[GU]{Aytekin \c{C}ibik\corref{cor1}}
\ead{abayram@gazi.edu.tr}

\author[UPitt]{Rui Fang}
\ead{ruf10@pitt.edu}

\author[UPitt]{William Layton}
\ead{wjl@pitt.edu}

\author[GaTech]{Farjana Siddiqua}
\ead{fsiddiqua3@gatech.edu}

\address[GU]{Department of Mathematics, Gazi University, Ankara, 06550, T\"{u}rkiye}
\address[UPitt]{Department of Mathematics, University of Pittsburgh, Pittsburgh, PA, 15260, USA}
\address[GaTech]{School of Mathematics, Georgia Institute of Technology, Atlanta, GA, 30309, USA}

\cortext[cor1]{Corresponding author. The work of William Layton and Rui Fang was partially supported by NSF grant DMS 2410893.}

    \begin{abstract}
Nudging is a data assimilation method amenable to both analysis and implementation. It also has the (reported) advantage of being insensitive to model errors compared to other assimilation methods. However, nudging behavior in the presence of model errors is little analyzed. This report gives an analysis of nudging to correct model errors. The analysis indicates that the error contribution due to the model error decays as the nudging parameter $\chi \to \infty$ like $\mathcal{O}(\chi^{-\frac{1}{2}})$, Theorem \ref{contis_model_error}. Numerical tests verify the predicted convergence rates and validate the nudging correction to model errors. 
    \end{abstract}

 
\end{frontmatter}
\section{Introduction}
Data assimilation aims to combine incomplete, noisy, and often time-delayed observations with a flow's incomplete and approximate dynamic laws to predict flow statistics more accurately than data or dynamic laws individually. Nudging, from \cite{luenberger1964observing}, is a data assimilation method with the (often reported but little analyzed) advantage of being insensitive to model errors.

This report begins the analysis of nudging to correct model errors. This analysis requires specifying a type of model error. We thus assume herein that the (true) flow data comes from a solution of 
\begin{equation}\label{nse_corio}
\frac{\partial u}{\partial t} + u\cdot \nabla u -\nu \Delta u +\nabla p + \omega R (u) = f, \nabla \cdot u=0.
\end{equation}
Here, $\omega R(u)$ represents the term omitted from the computational model, and $\omega$ is the scaling parameter of the term. For example, $\omega R(u)$ can represent a Coriolis or buoyant force, linked to an equation for temperature (see Section \ref{temperature_test} for details). The nudged approximation (with model error) takes sampled data, $I_H u$, chooses a nudging parameter $\chi\gg 1$, and solves a discretization of 
\begin{equation}\label{nudging}
\frac{\partial v}{\partial t} + v \cdot \nabla v -\nu \Delta v + \nabla \lambda - \chi I_H(u-v)=f, \ \nabla \cdot v =0.
\end{equation}
The formulation $(\ref{nse_corio})$, $(\ref{nudging})$ means that the model error \textit{residual}, $\|\omega R(v)\|$, can be computed. In Section \ref{error_wo_discre} the error $e= u(x,t)-v(x,t)$ is estimated under mild (abstract) conditions on $R(\cdot)$. In the cases easiest to analyze, we find that for small $H$ and (very) large $\chi$, the error due to inaccurate initial conditions decays exponentially (as in the case without model error). The error contribution due to model error decays as $\chi \to \infty$ like $\mathcal{O}(\chi^{-\frac{1}{2}})$, Theorem \ref{contis_model_error}. Section \ref{coriolis_case} makes a specific choice of $\omega R(\cdot)$, a Coriolis force term, and gives more details on the error behavior. In Section \ref{error_estimate_coriolis}, an analysis is given for the error between the space discretization of (\ref{nudging}) and the true $u(x,t)$. We find that the discrete solution $v^h$ of the nudging system converges to $u$ uniformly in time, and the model error term decreases as $\mathcal{O}{(\chi^{-\frac{1}{2}})}$ as $\chi \to \infty$.

Section \ref{temperature_test} presents extensive tests for a (slightly) more complex flow in which both (\ref{nse_corio}) and (\ref{nudging}) include temperature and buoyant forces, and (\ref{nse_corio}) includes a Coriolis force term. The tests show that nudging compensates for the absence of the Coriolis force in the solved model.

Problems (\ref{nse_corio}), (\ref{nudging}) are assumed to hold in $\Omega$, a bounded, convex, polyhedral domain in $2d$ or $3d$ with homogeneous Dirichlet (no-slip) boundary conditions. The initial conditions are specified for each problem.

\begin{equation}
u(x,0) = u_0(x), \text{ and } v(x,0) = v_0(x) \approx u_0(x), x \in \Omega.
\end{equation}
$I_H$ is an $L^2$ projection operator satisfying $\|(I-I_H)\phi\|\leq C_1 H \|\nabla \phi\|$ for all $\phi\in H^1(\Omega)$. The analysis is based on a condition on $H$ and $\chi$.  Equations (\ref{nse_corio}), (\ref{nudging}) are solved as a sequence of time intervals 
\begin{equation*}
0= t_0 <t_1 <t_2<\ldots,
\end{equation*}
where a $\chi = \chi_n$ value is selected by adapting each interval following \cite{CIBIK2025117526}. These values are chosen to satisfy
\begin{eqnarray}\label{condi_3d}
\chi _{n}-\omega q_{0}-\frac{2048}{19683}\nu ^{-3}\|\nabla v\|_{L^{\infty
}(t_{n},t_{n+1};L^{2}(\Omega ))}^{4} &\geq &\frac{1}{2}\alpha \chi _{n}>0, \\
\nu -2\left( C_{1}H\right)^{2}\chi_n  &>&0, \label{h_condi}
\end{eqnarray}
for some fixed $\alpha \in (0,1)$.

The conditions (\ref{condi_3d}), (\ref{h_condi}) above require stronger regularity than the analogous conditions (see (1.5), (1.6) in \cite{CIBIK2025117526} when model error is not considered). The way the stronger assumptions (\ref{condi_3d})-(\ref{h_condi}) are used in the analysis is detailed in the proof of Theorem \ref{contis_model_error} below. Reducing the assumed regularity from $\nabla v \in L^\infty(t_n, t_{n+1}; L^2)$ to $L^2(L^2)$ and $L^4(L^2)$ in $2d$ and $3d$, respectively, is an open problem. 

\subsection{Related Work}
Nudging has been used for data assimilation in many applications, e.g. \cite{CAO2022103659, FARHAT201559, Farhat_Martinez_2020, Jolly2019, Yuan_Pei, BOT15, Leo2022, FLT316, JMT17, ALT16, Farhat2016, GOT16, GLRVZ21, newey2024model}. The analysis of nudging without model errors is highly advanced in many papers, including the works of Azouani, Olson, Titi, and Edriss \cite{azouani2014continuous}; Biswas and Price \cite{Biswas_3D}; and Rebholz and Zerfas \cite{https://doi.org/10.1002/num.22751}. Extensive numerical experiments \cite{lei2015nudging, auroux2008nudging, dwight2010bayesian, chen2023data} demonstrate the robustness of nudging against model errors. The rigorous data assimilation analysis with model errors has received limited attention in the literature. Chen, Farhat, and Lunasin \cite{chen2023data} study assimilating Navier-Stokes equations (NSE) data into the Voight-$\alpha$-model (where the model error is the Voight term $\alpha^2\Delta v_t$). They prove in Theorem 4.8, p.9, that in $2d$ assimilation can decrease the model error from $\mathcal{O}(\alpha^2)$ to $\mathcal{O}(\alpha^4)$. Tests are given for a shell model of the Voight model. Building on this interesting work, we show model errors $\to 0$ in $2d$ and $3d$ for model errors with the structures (\ref{nse_corio}), (\ref{nudging}). These include physically motivated buoyancy and rotation \cite{lee2017solutions, kley1998treatment}.

\section{Notation and Preliminaries}
In this section, we introduce some of the notations and results used in this paper. We denote by $\|\cdot\|$ and $(\cdot,\cdot)$ the $L^2(\Omega)$ norm and inner product, respectively. We denote the $L^p(\Omega)$ norm by $\|\cdot\|_{L^p}$.
The solution spaces $X$ for the velocity and $Q$ for the pressure are defined as:
\begin{equation*}
\begin{aligned}
    &X:=(H_0^1(\Omega))^d=\{ v\in (L^2(\Omega))^d: \nabla  v\in (L^2(\Omega))^{d\cross d}\ \text{and}\  v=0\ \text{on}\ \partial\Omega\},\\
    &Q:=L^2_0(\Omega)=\{q\in L^2(\Omega): \int_\Omega q\ d x=0\}.
\end{aligned}
\end{equation*}
The divergence free subspace $V$ of $X$ is
\begin{equation*}
V=\{ v \in X: (q, \nabla \cdot v)=0, \forall q \in Q\}.
\end{equation*}
We denote Bochner Space \cite{adams2003sobolev} norm by $\|v\|_{L^p(0,T;X)}=\Bigg(\int_0^T\|v(\cdot,t)\|_{X}^p dt\Bigg)^{\frac{1}{p}}$, $p\in [1,\infty)$.

The space $H^{-1}(\Omega)$ denotes the dual space of bounded linear functionals defined on $X$ equipped with the norm:
\begin{equation*}
    \|f\|_{-1}=\sup_{0\neq  v\in X}\frac{(f, v)}{\|\nabla  v\|}.
\end{equation*}
The finite element method for this problem involves picking finite element spaces \cite{layton2008introduction}, $X^h\subset X$ and $Q^h\subset Q$. The discretely divergence-free subspace of $X^h$ is
\begin{equation*}
V^h :=\{v^h \in X^h: (\nabla \cdot v^h, q^h)=0, \ \forall q^h\in Q^h\}.
\end{equation*}

We assume that $(X^h, Q^h)$ are conforming. Motivated by \cite{girault1979finite} p. 91 error equation (1.29) and the discrete inf-sup condition (\ref{infsup}), we assume $V^h$ has approximation properties in  (\ref{prop}) similar to $X^h$, for $u\in (H^{m+1}(\Omega))^d\ \cap\ (H_0^1(\Omega))^d$ and $p\in H^{m}(\Omega)$,
\begin{equation}\label{prop}
\begin{aligned}
\inf_{v^h\in V^h}\{\|u-v^h\|+h\|\nabla(u-v^h)\|\}&\leq Ch^{m+1}|u|_{m+1},
\\\inf_{q^h\in Q^h}\|p-q^h\|&\leq Ch^{m}|p|_{m},
\end{aligned}
\end{equation}
\begin{equation}\label{infsup}
\inf_{q^h\in Q^h}\sup_{ v^h\in X^h}\frac{(q^h,\nabla\cdot  v^h)}{\|q^h\|\|\nabla  v^h\|}\geq \beta^h>0,
\end{equation}
where $\beta^h$ is bounded away from zero uniformly in $h$. We also use the inequalities derived by Ladyzhenskaya \cite
{ladyzhenskaya1969mathematical}.

\begin{lemma}(The Ladyzhenskaya Inequalities, Ladyzhenskaya \cite
{ladyzhenskaya1969mathematical}) For any vector function $u:{{\mathbb{R}}}%
^{d}\rightarrow {{\mathbb{R}}}^{d}$ with compact support and with the
indicated $L^{p}$ norms finite, 
\begin{eqnarray*}
\Vert u\Vert _{L^{4}({{\mathbb{R}}}^{2})} &\leq &2^{1/4}\Vert u\Vert^{1/2}\Vert \nabla u\Vert^{1/2},\ (d=2), \\
\newline
\Vert u\Vert _{L^{4}({{\mathbb{R}}}^{3})} &\leq &\frac{4}{3\sqrt{3}}\Vert
u\Vert ^{1/4}\Vert \nabla u\Vert ^{3/4},\ (d=3),\newline
\\
\Vert u\Vert _{L^{6}({{\mathbb{R}}}^{3})} &\leq &\frac{2}{\sqrt{3}}\Vert
\nabla u\Vert ,\ (d=3).
\end{eqnarray*}
\end{lemma} 

The standard explicitly skew-symmetrized trilinear form is denoted $$b^*(u,v,w)=\frac{1}{2}(u\cdot \nabla v,w)-\frac{1}{2}(u\cdot \nabla w, v).$$ Further, we have
\begin{equation}
    b^{\ast }(u,v,w) =(u\cdot \nabla v,w)+\frac{1}{2}(\left( \nabla \cdot
u\right) v,w),
\end{equation}
We define constants $M,K$ as follows%
\begin{eqnarray*}
K:=\sup_{u,v,w\in V}\frac{b^{\ast }(u,v,w)}{\Vert \nabla u\Vert \Vert
\nabla v\Vert \Vert \nabla w\Vert }<\infty , \
M :=\sup_{u,v,w\in X}\frac{b^{\ast }(u,v,w)}{\Vert \nabla u\Vert \Vert
\nabla v\Vert \Vert \nabla w\Vert }<\infty .
\end{eqnarray*}
\begin{lemma} \label{nonlinear_bound_John}
(Lemma 6.14, p. 311 - p. 312, \cite{john2016finite}) There is a $C_2<\infty$ with
\begin{equation}
\begin{split}
 b^*(u,v, w)
&\leq C_2\sqrt{\|u\|}\sqrt{\|\nabla u\|} \|\nabla v\|\|\nabla w\|,
 \end{split}
\end{equation}
for $u,v,w \in (H^1_0)^d$, $d=2$ or $3$.
\end{lemma}

\begin{lemma}(A discrete Gronwall lemma, Lemma 5.1, p. 369, \cite{heywood1990finite}) \label{discrete_gronwall} Let $\Delta t, B, a_n, b_n, c_n ,d_n$ be
nonnegative numbers such that for $l
\geq 1$:
\begin{equation*}
\begin{split}
a_l + \Delta t \sum_{n=0}^l b_n\leq \Delta t \sum_{n=0}^{l-1} d_n a_n + \Delta t \sum_{n=0}^{l}c_n + B, \text{ for } l\geq 0,
\end{split}
\end{equation*}
then for all $\Delta t>0$,
\begin{equation*}
\begin{split}
a_l + \Delta t \sum_{n=0}^l b_n\leq \exp(\Delta t \sum_{n=0}^{l-1} d_n) (\Delta t \sum_{n=0}^l c_n +B).
\end{split}
\end{equation*}
\end{lemma}
\section{Correction of model error}\label{error_wo_discre}
This section analyses the error between (1.1) and (1.2), exploiting the additive form of the model term omitted from the nudged system (1.2). The
analysis is under assumptions A1 - A4 below, that $I_{H}$ is a projection, $R(\cdot )$
is Lipschitz continuous, $v(x,t)$ is a strong solution, H is small enough and 
$\chi $\ is large enough. The error behavior under favorable assumptions, in this section, delineates reasonable expectations when discretization is added to the following sections.

A1 Assumption on $I_{H}$: $I_{H}$ is an $L^{2}$ projection satisfying, for
all $\phi \in \left( H_{0}^{1}(\Omega )\right) ^{d},$%
\begin{equation*}
\begin{split}
\|\phi -I_{H}\phi \|\leq C_{1}H\|\nabla \phi \|.
\end{split}
\end{equation*}

A2 Assumption that $v$ is a strong solution: $v(x,t)$ is a strong solution
of the nudged system with $\nabla v\in L^{\infty }(0,\infty ;L^{2}(\Omega
))$ with bound
\[
ess \sup_{0 < t < \infty} \frac{1}{|\Omega|} \int_{\Omega} |\nabla v|^2 \, dx \leq \left(\frac{U}{\lambda
_{T}}\right)^2 <\infty.
\]%
Here, $U, \lambda_{T}$ are constants with units Length/Time and Length,
respectively. $\lambda _{T}$\ is an analogue of the Taylor microscale.

A3 $R(\cdot )$ is Lipschitz continuous: For all $u,v \in \left(
L^2(\Omega)\right)^{d}$,
\[
\|R(u)-R(v)\|\leq q_{0}\|u-v\|.
\]

A4 $H$ small and $\chi$ large: In $d=3$ dimensions, for some $\alpha \in
(0,1)$ and on each $t_{n}\leq t\leq t_{n+1}$,
\begin{eqnarray*}
\chi _{n}-\omega q_{0}-\frac{2048}{19683}\nu ^{-3}\|\nabla v\|_{L^{\infty
}(t_{n},t_{n+1};L^{2}(\Omega ))}^{4} \geq \frac{1}{2} \alpha \chi _{n}>0, \\
\nu -2\left( C_{1}H\right)^{2}\chi_n >0.
\end{eqnarray*}

\begin{remark}
There is extensive analysis of the $2d$ case in the nudging literature. In $d=2
$ dimensions the $\chi $\ condition in A4 simplifies. In d=3 dimensions, A4
differs from the analogous condition without model error in 2 ways. First, $q_{0}$ is included in the condition, a minor change in the expected case of $q_{0}=O(1)$. Second, more regularity is assumed. Without model error the
factor $||\nabla v||_{L^{\infty }(t_{n},t_{n+1};L^{2}(\Omega ))}^{4}$\ is
replaced by $||\nabla v||_{L^{4}(t_{n},t_{n+1};L^{2}(\Omega ))}^{4}$. If the
latter is bounded, the former can be proven bounded by a bootstrap argument, but a Reynolds number (Re)-dependent constant appears.

The assumption A3 that the Nemitski operator $u\rightarrow R(u)$ is
globally Lipschitz on $L^{2}$ restricts the form of $R(\cdot )$ for large
arguments, p. 15 - p. 22 in \cite{ambrosetti1995primer}.  
A3 can be weakened and A2 can be assumed about $u$ rather
than $v$.

The parameter $\chi$ has units 1/Time. Let the large-scale turnover time be
denoted $T^{\ast}:=L/U$. The condition A4 on $\chi$\ scales like%
\[
\chi T^{\ast}\simeq\nu^{-3}L^6\left( \frac{U}
{\lambda_{T}}\right)  ^{4}\cdot\frac{L}{U}=Re^{3}\left(  \frac{L}{\lambda_{T}
}\right)^{4}.
\]
For fully developed turbulence, the Taylor microscale is roughly (Pope \cite{pope2000turbulent},
p. 200) $L/\lambda_{T}\simeq Re^{+1/2}.$ In this case $\chi T^{\ast}\simeq
Re^{5}$ is extremely large. For academic test problems, both the conditions on
$\chi$\ and the resulting scaling can be improved, \cite{CIBIK2025117526}.
\end{remark}

We prove the error has two components. One decays exponentially in time
while the other is $O(\chi^{-\frac{1}{2}})$.

\begin{theorem}\label{contis_model_error}
Suppose $d=3$ and A1, A2, A3, A4 hold. Then, on each $t_{n}\leq t\leq
t_{n+1},e(x,t)=u(x,t)-v(x,t)$ \ satisfies%
\begin{eqnarray*}
\|e(t_{n+1})\|^{2} &\leq &\exp \{-\alpha \chi
(t_{n+1}-t_{n})\}\|e(t_{n})\|^{2}+\chi ^{-1}\frac{\omega ^{2}}{\alpha }%
\overline{R}^{2}, \\
\text{where }\overline{R} &=&\min \left\{ \|R(u)\|_{L^{\infty
}(t_{n}, t_{n+1}; L^{2}(\Omega ))},\|R(v)\|_{L^{\infty
}(t_{n},t_{n+1}; L^{2}(\Omega))}\right\} .
\end{eqnarray*}
\end{theorem}

\begin{proof}
The proof begins with some standard steps for the NSE terms. Subtract (1.1)
and (1.2) and take the inner product with $e$, integrate by parts, and use
skew-symmetry of the nonlinearity. This gives
\begin{equation}\label{continous_error_inner}
\frac{1}{2}\frac{d}{dt}\|e\|^{2}+(e\cdot \nabla v,e)+\nu \|\nabla
e\|^{2}+\chi \|I_{H}e\|^{2}=\omega (R(u),e).
\end{equation}
By A1, $\chi \|I_{H}e\|^{2}\geq \chi \|e\|^{2}-\chi (C_{1}H)^{2}\|\nabla
e\|^{2}$. In $3d$ we use $(e\cdot \nabla v,e)\leq \|\nabla v\|\|e\|^2_{L^4}$, the Ladyzhenskaya inequalities and the arithmetic-geometric mean inequality (see \cite{CIBIK2025117526} Proposition 2.1 for
details) to arrive at
\[
(e\cdot \nabla v,e)\leq \frac{\nu }{2}\|\nabla e\|^{2}+\left( \frac{2048}{
19683}\nu ^{-3}\|\nabla v(t)\|^{4}\right) \|e\|^{2}.
\]
Combining the above, multiplying by 2, using A2, A4, and the arithmetic-geometric mean inequality again. By A4 $\nu - 2(C_1H)^2 \chi \geq 0$, it yields
\[
\frac{d}{dt}\|e\|^{2}+\left[ \nu -2\left( C_{1}H\right) ^{2}\chi \right]
\|\nabla e\|^{2}+\alpha \chi \|e\|^{2}\leq \chi ^{-1}\frac{\omega ^{2}}{%
\alpha }\|R(u)\|^{2}.
\]
Dropping the positive $\|\nabla e\|^{2}$\ term and using an integrating
factor gives $$\|e(t_{n+1})\|^{2}\leq \exp \{-\alpha \chi
(t_{n+1}-t_{n})\}\|e(t_{n})\|^{2}+\chi ^{-1}\frac{\omega ^{2}}{\alpha }%
\|R(u)\|^2_{L^{\infty }(t_{n},t_{n+1};L^{2}(\Omega ))}.$$
We now prove the same
estimate with $R(v$) on the right-hand side (RHS), allowing us to replace either with $\overline{R}$. In (\ref{continous_error_inner}) the RHS is estimated instead by
\begin{eqnarray*}
RHS &=&\omega (R(u),e)=\omega (R(u)-R(v),e)+\omega (R(v),e) \\
\text{(by A3)} &\leq &\omega q_{0}\|e\|^{2}+\frac{\chi \alpha }{2}\|e\|^{2}+\frac{\omega
^{2}}{2\alpha \chi } \|R(v)\|^{2}.
\end{eqnarray*}%
The remainder of the proof is the same.
\end{proof}

\section{The discrete problems}\label{coriolis_case}
The weak formulation of \eqref{nse_corio} is as follows: Find $( u,p)\in (X,Q)$ such that
\begin{align}
\label{nse_corio_v1}
    ( u_t, v)  +b^{*}(u,u,v)+\nu(\nabla u,\nabla  v)
    -( p,\nabla\cdot  v)+(\omega Ru,v)&=(f, v),\quad\forall\   v\in X,\\
    \label{nse_corio_v2}
    (\nabla \cdot  u,q)&=0, \quad\forall\ q\in Q.
\end{align}
The weak formulation of \eqref{nudging} is as follows: Find $(v,\lambda)\in (X,Q)$ such that $\forall\   w\in X$ and $\forall\ q\in Q$,
\begin{align}
\label{nudging_c_v1}
    ( v_t, v)  +b^{*}(v,v,w)+\nu(\nabla v,\nabla  w)
    -( \lambda,\nabla\cdot  w)+\chi(I_H(v-u),w)&=(f, w),\\
    \label{nudging_c_v2}
    (\nabla \cdot w,q)&=0.
\end{align}
The semi-discrete approximation of \eqref{nudging} is as follows.
Suppose $v^h(x,0)$ is an approximation of $v(x,0)$. The approximate velocity and pressure are maps
\begin{equation*}
     v^h: [0,T]\to X^h, \ \lambda^h:(0,T]\to Q^h 
\end{equation*}
\text{satisfying}, $\forall \ w^h\in X^h \ \text{and}\ \forall\  q^h\in Q^h$,
\begin{equation}\label{nudging_v1}
    ( v_t^h, w^h)  +b^{*}(v^h,v^h,w^h)+\nu(\nabla v^h,\nabla  w^h)
    -( \lambda^h,\nabla \cdot  w^h) +\chi (I_{H}(v^h-u),w^h)=(f, w^h),
    \end{equation}
    \begin{equation}\label{nudging_v2}
(\nabla \cdot v^h,q^h)=0.
    \end{equation}
The error in the semi-discrete approximation is analyzed in Section \ref{error_estimate_coriolis}. 

Section \ref{numerical_tests} presents numerical tests based on a linearly implicit Crank-Nicolson method, described next.  
We divide the time interval by $t_0=0<t_1<t_2<\cdots<t_N=T$. Let the time step and other quantities be denoted by 
\begin{align*}
    \text{time-step}=\Delta t,\ t_n=n\Delta t,\ f_n=f( x,t_n), \\
     v^{h}_n=\text{approximation to}\  v( x,t_n), \\
    \lambda^{h}_n=\text{approximation to} \ \lambda( x,t_n).
\end{align*}

We define $t_{n+1/2}=\frac{t_n+t_{n+1}}{2}$. Let $v^h_{n+1/2}$ be an approximation to $v(x,t_{n+1/2})$. For function $v$, we denote
\begin{equation*}
v^h_{n+1/2}=\frac{ v^h_n+ v^h_{n+1}}{2},\quad \Tilde{ v}^h_{n+1/2}=\frac{3 v^h_n- v^h_{n-1}}{2}.
\end{equation*}
To get a full discretization, we consider finite element spatial discretization and the second-order linearly implicit Crank-Nicolson scheme (also referred to as CNLE) for time discretization.
Given $( v^{h}_n,\lambda^{h}_n)\in (X^h,Q^h)$, for all $w^h\in X^h$ and $q^h\in Q^h$, find $( v^{h}_{n+1},\lambda^{h}_{n+1})\in (X^h,Q^h)$ satisfying
\begin{align}\label{nudging_fem_mid1}
  &\Big( \frac{ v^{h}_{n+1}- v^{h}_{n}}{\Delta t}, w^h\Big)+b^*( \widetilde{v}^{h}_{n+1/2}, v^{h}_{n+1/2}, w^h) 
    +\nu(\nabla  v^{h}_{n+1/2},\nabla  w^h)-(\lambda^{h}_{n+1/2},\nabla \cdot  w^h)\notag
    \\&+\chi (I_H(v^{h}_{n+1/2}-u(t_{n+1/2})),w^h)=(f_{n+1/2}, w^h),
    \\&
\label{nudging_fem_mid2}
    (\nabla \cdot v_{n+1/2}^h,q^h)=0.
    \end{align}  
\subsection{Stability when $R(u)$ is a Coriolis force}\label{stability}
Section \ref{error_wo_discre} gives an estimate for a general, Lipschitz $R(u)$. Here we analyze the error under the additional assumption that $\left(R(w), w\right)=0$ and $\|R(w)\|\leq \|w\|$ for any $w \in X$ (as for Coriolis forces). In this section, we prove the stability of the time-continuous nudged solution \( v^h \) in the finite element space. We also analyze the stability of the fully discrete nudged solution, which is computed using CNLE above.

First, since $(R(u),u)=0$, the standard energy estimate for the NSE (e.g. \cite{layton2008introduction}) also holds for \eqref{nse_corio}:
\begin{equation*}
\|u(T)\|^2 +\int_{0}^T \nu \|\nabla u\|^2 \leq \|u(0)\|^2 + \frac{1}{\nu}\int_{0}^T \|f\|^2_{-1}\, dt.
\end{equation*}

\begin{theorem} \label{thm:stability2}(Stability of $ v^h$) Assume A1 holds and $(R(w),w)=0$ for all $w\in X$. Then \eqref{nudging} is unconditionally stable. The semi-discrete solution $v^h$ of (\ref{nudging_v1}), (\ref{nudging_v2})satisfies
\begin{equation}
\begin{split}
\|v^h(T)\|^2 +\int_{0}^T (\nu \|\nabla  v^h\|^2+\chi\|I_H(v^h)\|^2+ \chi \|I_H(u-v^h)\|^2) \, dt \\
\leq \frac{1}{\nu} \int_{0}^T \|f\|^2_{-1}\, dt + \chi \int_{0}^T \|I_H(u)\|^2 \, dt.
\end{split}
\end{equation}
\end{theorem}
\begin{proof}
Take $w^h=v^h$ in \eqref{nudging_v1}, $q^h=\lambda ^h$ in \eqref{nudging_v2} and add them. We get the following energy equality,
\begin{equation*}
\frac{1}{2}\frac{d \|v^h\|^2}{dt} + \nu \|\nabla  v^h\|^2 + \chi (I_H(v^h),v^h)=(f,v^h)+\chi (I_H u,v^h).
\end{equation*}
For the nudging term, $\chi (I_H(v^h), v^h)=\chi\|I_H(v^h)\|^2$ because $I_H$ is an $L^2$ projection. Further, we have
\begin{equation*}
(f,v^h)\leq \frac{1}{2\nu}\|f\|^2_{-1} + \frac{\nu}{2}\|\nabla  v^h\|^2, \text{ and }
\end{equation*}
Next,
\begin{equation*}
\begin{gathered}
    \chi (I_H(u), v^h) = \chi (I_H(u), I_H(v^h)) = \frac{\chi}{2}\|I_H(u)\|^2 +\frac{\chi}{2}\|I_H(v^h)\|^2 -\frac{\chi}{2}\|I_H(u-v^h)\|^2.
\end{gathered}
\end{equation*}
Combining all, and multiplying by $2$, we have
\begin{equation*}
\frac{d\|v^h\|^2}{dt}+ \nu \|\nabla  v^h\|^2 + \chi \|I_H(v^h)\|^2 \leq \frac{1}{\nu} \|f\|^2_{-1} + \chi \|u\|^2.
\end{equation*}
Integrating from $0$ to $T$, we have the final result. 
\end{proof}

Now, we prove the stability of the CNLE time-discretization.
\begin{theorem}\label{thm:stability4} (CNLE stability) Assume A1 holds and $(R(w),w)=0$ for all $w\in X$. Then
\eqref{nudging_fem_mid1} is unconditionally energy stable. For any $N\geq 1$,
\begin{equation*}
\begin{gathered}
\|v^{h}_{N}\|^2  +\nu \Delta t \sum_{n=0}^{N-1}\|\nabla  v^{h}_{n+1/2}\|^2 +
\Delta t \chi \sum_{n=0}^{N-1} \|I_H(v^{h}_{n+1/2})\|^2+\Delta t \chi \sum_{n=0}^{N-1}\|I_H(u(t_{n+\frac{1}{2}})-v^h_{n+\frac{1}{2}})\|^2 \\
\leq \|v^{h}_{0}\|^2 +\frac{\Delta t}{\nu} \sum_{n=0}^{N-1}\|f_{n+1/2}\|^2_{-1} + 
\Delta t\chi \sum_{n=0}^{N-1}\|u(t_{n+1/2})\|^2.
\end{gathered}
\end{equation*}
\end{theorem}
\begin{proof}
Take $w^h=v^h_{n+1/2}$ in \eqref{nudging_fem_mid1}, $q^h=\lambda^{h,n+1/2}$ in \eqref{nudging_fem_mid2} and add them. $b^*(\Tilde{ v}^h_{n+\frac{1}{2}},  v^h_{n+\frac{1}{2}}, v^h_{n+\frac{1}{2}})=0$. Using other techniques used in \Cref{thm:stability2}, we get the following inequality,
\begin{equation*}
\begin{gathered}
\frac{1}{2}\|v^{h}_{n+1}\|^2 -\frac{1}{2}\|v^{h}_{n}\|^2 + \Delta t \nu \|\nabla  v^{h}_{n+1/2}\|^2+
\Delta t \chi\|I_H(v^{h}_{n+1/2})\|^2
\\\leq \frac{\Delta t}{2\nu} \|f_{n+1/2}\|^2_{-1}+\frac{\nu \Delta t}{2}\|\nabla  v^{h}_{n+1/2}\|^2 
+\frac{\Delta t \chi}{2} \|I_H(u(t_{n+1/2}))\|^2 \\
+ \frac{\chi \Delta t}{2}\|I_H(v^{h}_{n+1/2})\|^2-\frac{\Delta t \chi }{2}\|I_H(u(t_{n+\frac{1}{2}})-v^h_{n+\frac{1}{2}})\|^2.
\end{gathered}
\end{equation*}
Multiply by $2$, we have
\begin{equation*}
\begin{gathered}
\|v^{h}_{n+1}\|^2-\|v^{h}_{n}\|^2 + \Delta t \nu \|\nabla  v^{h}_{n+1/2}\|^2+
\Delta t \chi\|I_H(v^{h}_{n+1/2})\|^2
+\Delta t \chi \|I_H(u(t_{n+\frac{1}{2}})-v^h_{n+\frac{1}{2}})\|^2\\\leq \frac{\Delta t}{\nu} \|f_{n+1/2}\|^2_{-1}
+\Delta t \chi \|u(t_{n+1/2})\|^2.
\end{gathered}
\end{equation*}
Summing from $n=0$ to $n=N-1$, we have the final result. 
\end{proof}

\section{Error estimates for spatial discretization when $R(u)$ is a Coriolis force}\label{error_estimate_coriolis}
In this section, we analyze the error estimates under the additional assumption that $\left(R(w), w\right)=0$ and $\|R(w)\|\leq \|w\|$ for any $w\in X$ (as in the setting considered in Section \ref{stability}). We present the error estimates for the semi-discrete case. Define the error $e:= u-v^h$. We decompose the error into two pieces: $e= \eta -\phi^h$, where $\eta = u- \Tilde{u}$, and $\phi = v^h-\Tilde{u} \in V^h$, $\Tilde{u}\in V^h$. We assume a discrete condition that is analogous to A4: for some $\alpha \in (0,1)$, on each $t_n\leq t\leq t_{n+1}$,
\begin{equation}\label{weaker_A4}
\begin{split}
\chi_n -C_2^4\left(\frac{9261}{8}\right) \nu^{-3} \|\nabla v^h\|^4_{L^\infty(0,\infty; L^2(\Omega))}&\geq \alpha \chi_n>0, \\
\nu-C_1^2 H^2 \chi_n &>0.
\end{split}
\end{equation}

Herein, we bound the error in terms of $u$, which does not rely on A2. Hence no $q_0$ appears in (\ref{weaker_A4}). The constant $C_2^4 \frac{9261}{8}$ is tailored to the discrete case with the skew-symmetrized trilinear term. In the next proof, we fix $\chi$ to shorten the analysis. For variable $\chi = \chi_n$, perform the proof steps on each $t_n\leq t\leq t_{n+1}$, then sum over (see \cite{CIBIK2025117526} for this plan executed). 
\begin{theorem} \label{theorem_error} Let $\widetilde{u}\in V^h$ be an approximation of $u$. Assume A1, (\ref{weaker_A4}), and $\|\nabla u\| \in L^4(0,T)$ hold, we have
\begin{equation}\label{error-semi}
\begin{gathered}
\|e(t)\|^2 +\frac{\nu}{2}\int_{0}^t \exp{-\frac{\alpha}{2} \chi (t-t')} \|\nabla e\|^2\, dt' \leq \exp{\frac{-\alpha}{2} \chi t}\|e(0)\|^2 \\
+ \frac{2}{\alpha} \max_{0\leq t\leq T}\|I_H(u-\Tilde{u})\|^2+2 \omega^2\left(\frac{1}{\chi}+ \frac{7C_1^2}{2\nu} H^2 \right) \|u\|^2_{L^2(0,T;L^2)}\\
+\frac{7}{\nu}K^2\left(\|\nabla v^h\|^2_{L^4(0,T;L^2)}+ \|\nabla u\|^2_{L^4(0,T;L^2)}\right) \|\nabla (u-\Tilde{u})\|^2_{L^4(0,T;L^2)}+\frac{7}{\nu}\int_{0}^t \|p-q^h\|^2 
\, dt'\\
+\frac{7}{\nu}\int_{0}^t \|(u-\Tilde{u})_t\|^2_{-1}\, dt'+ \| (u-\Tilde{u})(t)\|^2+ 8\nu \|\nabla (u-\Tilde{u})\|^2_{L^2(0,T;L^2)}
\end{gathered}
\end{equation}
\end{theorem}
\begin{proof}
Take the difference of equation (\ref{nse_corio_v1}) and equation (\ref{nudging_v1}), and take $w^h=\phi^h$ it yields
\begin{equation}
\begin{gathered}
\frac{1}{2}\frac{d \| \phi^h\|^2}{d t} + \nu \|\nabla \phi^h\|^2 +\chi (I_H(v^h-u), \phi^h) -\omega (R u,\phi^h)+(p,\nabla \cdot \phi^h) \\
\leq b^*(u, u,\phi^h) -b^*(v^h,v^h,\phi^h) +(\eta_t, \phi^h) + \nu (\eta, \phi^h).
\end{gathered}
\end{equation}
Consider the nudging term $\chi(I_H(v^h-u), \phi^h)$:
\begin{equation}
\begin{split}
\chi(I_H(v^h-u), \phi^h)
&=\chi (I_H(\phi^h), \phi^h) -\chi (I_H(\eta),\phi^h)\\
&=\chi \|I_H(\phi^h)\|^2 -\chi (I_H(\eta),I_H(\phi^h))\\
&\geq \chi \|I_H(\phi^h)\|^2- \frac{\chi}{2}\|I_H(\eta)\|^2 - \frac{\chi}{2}\|I_H(\phi^h)\|^2\\
&\geq \frac{\chi}{2}\|I_H(\phi^h)\|^2- \frac{\chi}{2}\|I_H(\eta)\|^2.
\end{split}
\end{equation}
The model error term, $\omega(R u,\phi^h)$ satisfies
\begin{equation*}
\omega(R u,\phi^h)\leq\omega \|R(u)\|\|\phi^h\|\leq \omega \|u\|\|\phi^h\|.
\end{equation*}
We use the triangle inequality and A1:
\begin{equation*}
\|\phi^h\|-\|I_H(\phi^h)\|\leq \|I_H(\phi^h)-\phi^h\|\leq C_1H \|\nabla \phi^h\|,
\end{equation*}
thus 
\begin{equation*}
\|\phi^h\|\leq \|I_H(\phi^h)\|+ C_1H\|\nabla \phi^h\|.
\end{equation*}
Hence we have
\begin{equation*}
\begin{gathered}
\omega(R u,\phi^h)\leq \omega \|u\|\|\phi^h\|\leq \omega \|u\|\|I_H(\phi^h)\|+ \omega C_1H \|u\|\|\nabla \phi^h\|\\
\leq \omega^2 \left(\frac{1}{\chi}+\frac{7C_1^2}{2\nu}H^2  \right) \|u\|^2 + \frac{\chi}{4}\|I_H(\phi^h)\|^2 + \frac{\nu}{14}\|\nabla \phi^h\|^2.
\end{gathered}
\end{equation*}
The NSE terms are treated as usual by
\begin{equation}
(\eta_t, \phi^h)\leq \|\eta_t\|_{-1}\|\nabla \phi^h\|\leq \frac{\nu}{14}\|\nabla \phi^h\|^2 +\frac{7}{2\nu}\|\eta_t\|^2_{-1}.
\end{equation}
\begin{equation*}
\nu (\nabla \eta, \nabla \phi^h)\leq \nu\|\nabla \eta\|\|\nabla \phi^h\|\leq \frac{\nu}{14} \|\nabla \phi^h\|^2 + \frac{7}{2}\nu \|\nabla \eta\|^2, \text{ and }
\end{equation*}
\begin{equation*}
\begin{gathered}
(p,\nabla \cdot \phi^{h})=(p-q^h, \nabla \cdot \phi^h)\\
\leq \|p-q^h\|\|\nabla \cdot \phi^h\|\leq \frac{\nu}{14}\|\nabla \phi^h\|^2 +\frac{7}{2\nu}\|p-q^h\|^2.
\end{gathered}
\end{equation*}
The nonlinear terms $b^*(\cdot, \cdot, \cdot)$ satisfies
\begin{equation*}
    \begin{gathered}
    b^*(u,u,\phi^h)- b^*(v^h,v^h,\phi^h)\\
        =b^*(u,u,\phi^h)+b^*(\phi^h-v^h,v^h,\phi^h)-b^*(\phi^h, v^h,\phi^h)\\
        =b^*(u,u,\phi^h) + b^*(\eta -u, v^h,\phi^h) -b^*(\phi^h, v^h,\phi^h)\\
        =b^*(u, \eta,\phi^h) + b^*(\eta, v^h,\phi^h) -b^*(\phi^h, v^h, \phi^h).
    \end{gathered}
\end{equation*}

Thus, term by term, 
\begin{equation}
\begin{split}
    b^*(u,\eta,\phi^h)
    \leq K \|\nabla \eta\| \|\nabla u\|\|\nabla \phi^h\|\leq \frac{\nu}{14}\| \nabla \phi^h\|^2+ \frac{7}{2\nu}K^2\|\nabla u\|^2 \|\nabla \eta\|^2.
\end{split}
\end{equation}
\begin{equation}
\begin{gathered}
    b^*(\eta, v^h,\phi^h)\leq  K\|\nabla \eta\| \|\nabla v^h\|\|\nabla \phi^h\|
    \leq \frac{\nu}{14} \|\nabla \phi^h\|^2 +\frac{7}{2\nu}K^2 \|\nabla v^h\|^2\|\nabla \eta\|^2
\end{gathered}
\end{equation}
By Lemma \ref{nonlinear_bound_John}, we have
\begin{equation}
\begin{gathered}
    b^*(\phi^h, v^h,\phi^h)\leq C_2\|\phi^h\|^{1/2}\|\nabla \phi^h\|^{3/2}\|\nabla v^h\|\\
    \leq \frac{\nu}{14}\|\nabla \phi^h\|^2+C_2^4\frac{9261}{32}\nu^{-3}\|\nabla v^h\|^4\|\phi^h\|^2.
\end{gathered}
\end{equation}

Combining the above, multiplying by $2$ and simplifying gives
\begin{equation}\label{error_put_together}
\begin{gathered}
\frac{d\|\phi^h\|^2}{dt} +\nu\|\nabla \phi^h\|^2  + \frac{\chi}{2}\|I_H(\phi^h)\|^2 \\
\leq \chi\|I_H(\eta)\|^2 + 2 \omega^2\left(\frac{1}{\chi}+ \frac{7C_1^2}{2\nu} H^2 \right) \|u\|^2\\
+ \frac{7}{\nu}\|\eta_t\|^2_{-1} + \left(7\nu+ \frac{7}{\nu}K^2 \|\nabla v^h\|^2+ \frac{7}{\nu}K^2\|\nabla u\|^2 \right) \|\nabla \eta\|^2 \\
+C_2^4\frac{9261}{16}\nu^{-3} \|\nabla v^h\|^4 \|\phi^h\|^2 + \frac{7}{\nu}\|p-q^h\|^2.
\end{gathered}
\end{equation}
As before, using A1 gives
\begin{equation*}
\begin{gathered}
\|I_H(\phi^h)\|^2 = \|\phi^h\|^2 - \|(I-I_H)(\phi^h)\|^2
\geq \|\phi^h\|^2 - C_1^2H^2 \|\nabla \phi^h\|^2.
\end{gathered}
\end{equation*}
Thus, we have
\begin{equation*}
\begin{gathered}
\frac{d \|\phi^h\|^2}{dt} + \frac{1}{2}\left(\nu-C_1^2 \chi H^2\right)\|\nabla \phi^h\|^2+\frac{\nu}{2}\|\nabla \phi^h\|^2 + \left(\frac{\chi}{2} - C_2^4\frac{9261}{16}\nu^{-3}\|\nabla v^h\|^4 \right)  \|\phi^h\|^2 \\
 \leq \chi\|I_H(\eta)\|^2 +2 \omega^2\left(\frac{1}{\chi}+ \frac{7C_1^2}{2\nu} H^2 \right) \|u\|^2
+ \frac{7}{\nu}\|\eta_t\|^2_{-1}\\ + \left(7\nu + \frac{7}{\nu}K^2\|\nabla v^h\|^2+\frac{7}{\nu}K^2\|\nabla u\|^2 \right) \|\nabla \eta\|^2 
+\frac{7}{\nu}\|p-q^h\|^2.
\end{gathered}
\end{equation*}
By (\ref{weaker_A4}), 
\begin{equation*}
\begin{gathered}
\frac{d \|\phi^h\|^2}{dt} +\frac{\nu}{2}\|\nabla \phi^h\|^2 + \frac{1}{2} \alpha \chi \|\phi^h\|^2 \\
\leq \chi\|I_H(\eta)\|^2 +2 \omega^2\left(\frac{1}{\chi}+ \frac{7C_1^2}{2\nu} H^2 \right) \|u\|^2
+ \frac{7}{\nu}\|\eta_t\|^2_{-1}\\ + \left(7\nu+\frac{7}{\nu}K^2\|\nabla v^h\|^2+ \frac{7}{\nu}K^2\|\nabla u\|^2 \right) \|\nabla \eta\|^2 
+\frac{7}{\nu}\|p-q^h\|^2.
\end{gathered}
\end{equation*}
Multiply both sides by the integration factor $\exp{\frac{\alpha}{2} \chi t}$. We integrate from $0$ to $t$, and multiply $\exp\{-\frac{\alpha}{2}\chi t\}$ to both sides:
\begin{equation*}
\begin{gathered}
\|\phi^h(t)\|^2 + \frac{\nu}{2}\int_{0}^t \exp{-\frac{\alpha}{2} \chi (t-t')}\|\nabla \phi^h\|^2\, dt'
\leq \exp{\frac{-\alpha}{2} \chi t} \|\phi^h(0)\|^2 \\
+ \int_{0}^t \exp{-\frac{\alpha}{2} \chi (t-t')} \bigg\{\chi\|I_H(\eta)\|^2 +2 \omega^2\left(\frac{1}{\chi}+ \frac{7C_1^2}{2\nu} H^2 \right) \|u\|^2
+ \frac{7}{\nu}\|\eta_t\|^2_{-1}\\
+ \left(7\nu+\frac{7}{\nu}K^2\|\nabla v^h\|^2 +\frac{7}{\nu}K^2\|\nabla u\|^2 \right) \|\nabla \eta\|^2 
+\frac{7}{\nu}\|p-q^h\|^2 
\bigg\}\, dt'.
\end{gathered}
\end{equation*}
We bound $\int_{0}^t\exp{-\frac{\alpha}{2} \chi (t-t')} \chi \|I_H(\eta)\|^2\, dt'$:
\begin{equation*}
\begin{split}
&\int_{0}^t\exp{-\frac{\alpha}{2} \chi (t-t')} \chi \|I_H(\eta)\|^2\, dt'\\
&\leq \max_{0\leq t\leq T} \|I_H(\eta)\|^2 \chi \exp{-\frac{\alpha}{2} \chi t} \int_{0}^t \exp{\frac{\alpha}{2} \chi t'} \, dt'\\
&=\max_{0\leq t\leq T} \|I_H(\eta)\|^2 \chi \exp{-\frac{\alpha}{2} \chi t} \frac{\exp{\frac{\alpha}{2}\chi t}-1}{\alpha \chi/2}\\
&=\max_{0\leq t\leq T} \|I_H(\eta)\|^2 \frac{2}{\alpha} (1-\exp{-\alpha/2\chi})
\leq \frac{2}{\alpha}\max_{0\leq t\leq T} \|I_H(\eta)\|^2.
\end{split}
\end{equation*}
Since  $t>t'$, and $t \leq T$, $\exp{-\frac{\alpha}{2} \chi(t-t')}\leq 1$. We have
\begin{equation*}
\begin{gathered}
\|\phi^h(t)\|^2 + \frac{\nu}{2}\int_{0}^t \exp{-\frac{\alpha}{2} \chi (t-t')}\|\nabla \phi^h\|^2\, dt'
\leq \exp{\frac{-\alpha}{2} \chi t} \|\phi^h(0)\|^2 \\
+ \frac{2}{\alpha} \max_{0\leq t\leq T}\|I_H(\eta)\|^2+ \int_{0}^t 2 \omega^2\left(\frac{1}{\chi}+ \frac{7C_1^2}{2\nu} H^2 \right) \|u\|^2
+ \frac{7}{\nu}\|\eta_t\|^2_{-1}\\
+ \left(7\nu+\frac{7}{\nu}K^2\|\nabla v^h\|^2+ \frac{7}{\nu}K^2\|\nabla u\|^2 \right) \|\nabla \eta\|^2 
+\frac{7}{\nu}\|p-q^h\|^2
\bigg\}\, dt'.
\end{gathered}
\end{equation*}
We simplify by the following estimates.
\begin{equation}\label{estimate_1}
    \begin{gathered}
        \int_{0}^t \|\nabla u\|^2 \|\nabla \eta\|^2\, dt' \leq \left(\int_{0}^t  \|\nabla u\|^4 \, dt' \right)^{1/2} \left(\int_{0}^t \|\nabla \eta\|^4\, dt'\right)^{1/2}\\
        \leq \|\nabla u\|^2_{L^4(0,T;L^2)}\|\nabla \eta\|^2_{L^4(0,T;L^2)}.
    \end{gathered}
\end{equation}

\begin{equation}\label{estimate_2}
    \begin{gathered}
        \int_{0}^t \|\nabla v^h\|^2 \|\nabla \eta\|^2\, dt' \leq \left(\int_{0}^t  \|\nabla v^h\|^4 \, dt' \right)^{1/2} \left(\int_{0}^t \|\nabla \eta\|^4\, dt'\right)^{1/2}\\
        \leq \|\nabla v^h\|^2_{L^4(0,T;L^2)}\|\nabla \eta\|^2_{L^4(0,T;L^2)}.
    \end{gathered}
\end{equation}

Hence,
\begin{equation}\label{estimate_3}
\begin{gathered}
\|\phi^h(t)\|^2 + \frac{\nu}{2}\int_{0}^t \exp{-\frac{\alpha}{2} \chi (t-t')}\|\nabla \phi^h\|^2\, dt'
\leq \exp{\frac{-\alpha}{2} \chi t} \|\phi^h(0)\|^2 \\
+ \frac{2}{\alpha} \max_{0\leq t\leq T}\|I_H(\eta)\|^2+2 \omega^2\left(\frac{1}{\chi}+ \frac{7C_1^2}{2\nu} H^2 \right) \|u\|^2_{L^2(0,T;L^2)}\\
+ \frac{7}{\nu}K^2  \left( \|\nabla v^h\|^2_{L^4(0,T;L^2)}+\|\nabla u\|^2_{L^4(0,T;L^2)}\right) \|\nabla \eta\|^2_{L^4(0,T;L^2)}\\
+\frac{7}{\nu}\int_{0}^t \|p-q^h\|^2 
\, dt'+\frac{7}{\nu}\int_{0}^t \|\eta_t\|^2_{-1}\, dt'
+ 7\nu \|\nabla \eta\|^2_{L^2(0,T;L^2)}.
\end{gathered}
\end{equation}
The triangle inequality, $\|e\|\leq \|\eta\|+ \|\phi^h\|$ gives
\begin{equation*}
\begin{gathered}
\|e(t)\|^2 +\frac{\nu}{2}\int_{0}^t \exp{-\frac{\alpha}{2} \chi (t-t')}\|\nabla e(t')\|^2\, dt' \leq \|\eta(t)\|^2 \\
+ \nu \int_{0}^t \|\nabla \eta\|^2\, dt'
+\|\phi^h(t)\|^2 +\frac{\nu}{2}\int_{0}^t \exp{-\frac{\alpha}{2} \chi (t-t')}\|\nabla \phi^h\|^2\, dt'.
\end{gathered}
\end{equation*}
Thus, we have the final result of (\ref{error-semi}).
\end{proof}
\begin{corollary}
Let $\widetilde{u}\in V^h$ be an approximation of $u$. Assume A1, (\ref{weaker_A4}), $\|\nabla u\| \in L^4(0,T)$, the approximation properties in (\ref{prop}), and for any $w\in H^1(\Omega)^d$,
\begin{equation}\label{property_I_H}
    \|I_H(w)\|\leq C_3 \|w\|.
\end{equation}
We have, for any $0<t\leq T$,
\begin{equation}\label{error_approximation}
\begin{gathered}
\|e(t)\|^2 +\frac{\nu}{2}\int_{0}^t \exp{-\frac{\alpha}{2} \chi (t-t')} \|\nabla e\|^2\, dt' \leq \exp{\frac{-\alpha}{2} \chi t}\|e(0)\|^2 \\
+ C_3 C\frac{2}{\alpha} h^{2m+2}\max_{0\leq t\leq T}|u(t)|^2_{m+1}+2 \omega^2\left(\frac{1}{\chi}+ \frac{7C_1^2}{2\nu} H^2 \right) \|u\|^2_{L^2(0,T;L^2)}\\
+\frac{7}{\nu}K^2\left(\|\nabla v^h\|^2_{L^4(0,T;L^2)}+ \|\nabla u\|^2_{L^4(0,T;L^2)}\right) Ch^{2m}\sqrt{\int_{0}^T |u|^4_{m+1}\, dt}
+\frac{7}{\nu}C h^{2m}\int_{0}^T |p|^2_m\, dt\\
+Ch^{2m+4}\frac{7}{\nu}\int_{0}^t |u_t|^2_{m+1}\, dt'+ Ch^{2m+2}|u|^2_{m+1} + 8\nu Ch^{2m}\int_{0}^T |u|^2_{m+1}\, dt.
\end{gathered}
\end{equation}
\end{corollary}
\begin{proof}
We now bound the terms on the RHS of (\ref{error-semi}) using the approximation properties in (\ref{prop}).

From (\ref{property_I_H}) and (\ref{prop}), we get:
\begin{equation}
\|I_H(u - \Tilde{u})\| \leq C_3 \|u - \Tilde{u}\| \leq C C_3 h^{m+1} |u|_{m+1}.
\end{equation}

For the $\|\nabla(u - \Tilde{u})\|^2_{L^4(0,T;L^2)}$ term, we have:
\begin{equation}
\|\nabla(u - \Tilde{u})\|^2_{L^4(0,T;L^2)} = \sqrt{ \int_{0}^T \|\nabla(u - \Tilde{u})\|^4\, dt } \leq C h^{2m} \sqrt{ \int_{0}^T |u|^4_{m+1}\, dt }.
\end{equation}

For the pressure term:
\begin{equation}
\int_{0}^T \|p - q^h\|^2\, dt \leq C h^{2m} \int_{0}^T |p|^2_m\, dt.
\end{equation}

For the $L^2$ norm and $H^1$-seminorm of $u - \Tilde{u}$:
\begin{align}
\|u - \Tilde{u}\| &\leq C h^{m+1} |u|_{m+1}, \\
\|\nabla(u - \Tilde{u})\|^2_{L^2(0,T;L^2)} &= \int_{0}^T \|\nabla(u - \Tilde{u})\|^2\, dt \leq C h^{2m} \int_{0}^T |u|^2_{m+1}\, dt.
\end{align}

Finally, for the dual norm of the time derivative:
\begin{equation}
\begin{aligned}
\|(u - \Tilde{u})_t\|_{-1} &= \sup_{0 \neq v \in X} \frac{((u - \Tilde{u})_t, v)}{\|\nabla v\|} = \sup_{0 \neq v \in X} \inf_{v^h \in X^h} \frac{((u - \Tilde{u})_t, v - v^h)}{\|\nabla v\|} \\
&\leq \|(u - \Tilde{u})_t\| \sup_{0 \neq v \in X} \inf_{v^h \in X^h} \frac{\|v - v^h\|}{\|\nabla v\|}\leq \|(u - \Tilde{u})_t\| \sup_{0 \neq v \in X} \inf_{v^h \in V^h} \frac{\|v - v^h\|}{\|\nabla v\|}.
\end{aligned}
\end{equation}

Using the approximation property in (\ref{prop}) with $m = 0$, we get:
\begin{equation}
\inf_{v^h \in V^h} \|v - v^h\| \leq C h \|\nabla v\|,
\end{equation}
and therefore,
\begin{equation}
\|(u - \Tilde{u})_t\|_{-1} \leq C h \|(u - \Tilde{u})_t\| \leq C h^{m+2} |u_t|_{m+1}.
\end{equation}
\end{proof}

\section{Numerical Experiments}\label{numerical_tests}
Section \ref{test: convergence_rate} verifies the convergence rate with respect to the timestep predicted in Theorem \ref{theorem_error}. In Section \ref{temperature_test}, we test nudging's correction of model error for a coupled NSE system for natural convection in a double-pane window problem. The data includes a Coriolis term while the nudged system does not.
\subsection{Convergence rate study}\label{test: convergence_rate}
We first check if the expected convergence rates in $\Delta t$ are observed. While Section \ref{error_estimate_coriolis} did not analyze the discrete time model, this step checks code correctness. This step uses a manufactured true solution. The results (below) show that added terms do not deteriorate the expected convergence orders. To do so, we pick 
\begin{equation*}
u(x,y,t)=e^{t}(\cos y,\sin x), \ \text{ and }
	p(x,y,t)=(x-y)(1+t)
\end{equation*}
as true solutions and insert these values into (\ref{nse_corio}) to calculate the body force term. We choose $\nu=1$, $\omega=1$, and $\chi=100$ for simplicity. The true solution values are projected into a $P0$ finite element space, which will act as a coarse mesh to get samples. For the approximate solutions, we consider a finer mesh with a barycenter refined triangulation along with Scott-Vogelius finite elements. Other tests (not reported herein) with a Taylor-Hood pair with skew-symmetric trilinear forms on a usual Delaunay triangulation yielded analogous results. The spatial mesh width $h$ is taken as $1/32$, which results in $43266$ total number of degrees of freedom (dofs), and the time interval is $[0,2]$. The CNLE time discretization (\ref{nudging_fem_mid1})-(\ref{nudging_fem_mid2}) is used here. Thus, a second-order in time convergence is expected. We use the $L^2$ norm to measure the error, and the results are given in Table \ref{tab:error-main}.
\begin{table}[H]
	\centering
		\begin{tabular}{||c|c|c|c||}
			\hline
			$\Delta t$ & error & rate \\
			\hline
			1       & 0.0021   & -       \\
			1/2     & 0.00052   & 2.01     \\
			1/4     & 0.00015  & 1.80    \\
			1/8     & 4.2e-5   & 1.90    \\
			1/16    & 1.1e-5   & 1.93    \\
            1/32    & 3.0e-6   & 1.87    \\
			\hline
		\end{tabular}
	\caption{Errors and rates of convergence. A second-order error in time is obtained as expected.}
	\label{tab:error-main}
\end{table}

\subsection{Double Pane Window Problem}\label{temperature_test}
In this test, we nudge velocity alone. This choice of assimilating only velocity is inspired by the important results in Farhat, Jolly, and Titi \cite{FARHAT201559}; and Farhat, Glatt-Holtz, Martinez, McQuarrie, and Whitehead \cite{Farhat_Martinez_2020}.
In this part, we consider a coupled NSE system instead of the pure NSE, a natural convection problem with rotation. The double pane window problem is a classical benchmark test of  \cite{dvd83}, which is used widely to test the performance of natural convection codes. The nudged system of equations solved here is:
\begin{eqnarray}\label{nc}
\frac{\partial w}{\partial t} + w \cdot \nabla w -Pr \Delta w + \nabla \lambda + \chi I_H(w-u)=Pr Ra T\textbf{g} + f,
\\
\nabla \cdot w =0,
\\
\frac{\partial T}{\partial t} + w \cdot \nabla T - \Delta T = \gamma.
\end{eqnarray}
Here $Pr$ is the Prandtl number, $Ra$ is the Rayleigh number, $\textbf{g}$ is the unit vector in gravitational direction, $f$ is the body force, and $\gamma$ is the heat source.

For brevity, we do not give the discretization details for this system, since it is very similar to the one described for NSE. The domain is a unit square with no-slip boundary conditions, horizontal walls are kept adiabatic, and vertical walls are kept at different constant temperatures. The fluid is assumed to be air, which has $Pr=0.71$. We use the second-order backward differentiation formula (BDF2) temporal discretization with linear extrapolation in the nonlinear term analogous to (4.7), and a Scott-Vogelius finite element pair with a barycenter refined mesh. The test runs until it reaches a steady state with $\Delta t =10^{-3}$. 

We will obtain samplings from a direct numerical simulation (DNS) \textbf{ with a Coriolis term} and nudge the natural convection system without a Coriolis term, (6.1) above. We examine whether nudging takes account of the model error due to the omitted Coriolis term. We consider two different cases, $Ra=10^4$ and $Ra=10^5$.

\subsubsection{\textbf{Case I: $Ra=10^4$}}
\begin{figure}[H]
\begin{subfigure}{0.3\textwidth}
\includegraphics[width=0.9\linewidth]{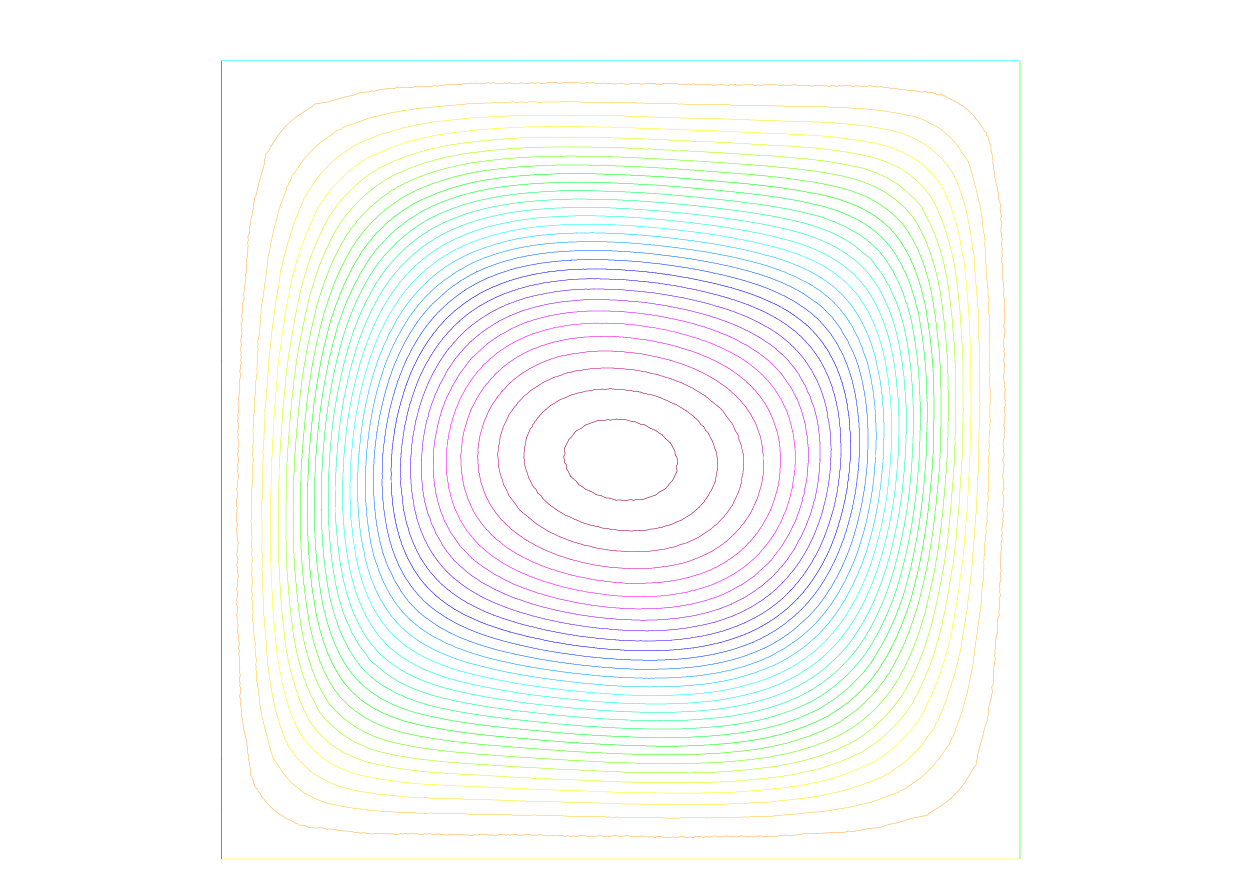} 
\caption{Streamlines}
\end{subfigure}
\begin{subfigure}{0.3\textwidth}
\includegraphics[width=0.9\linewidth]{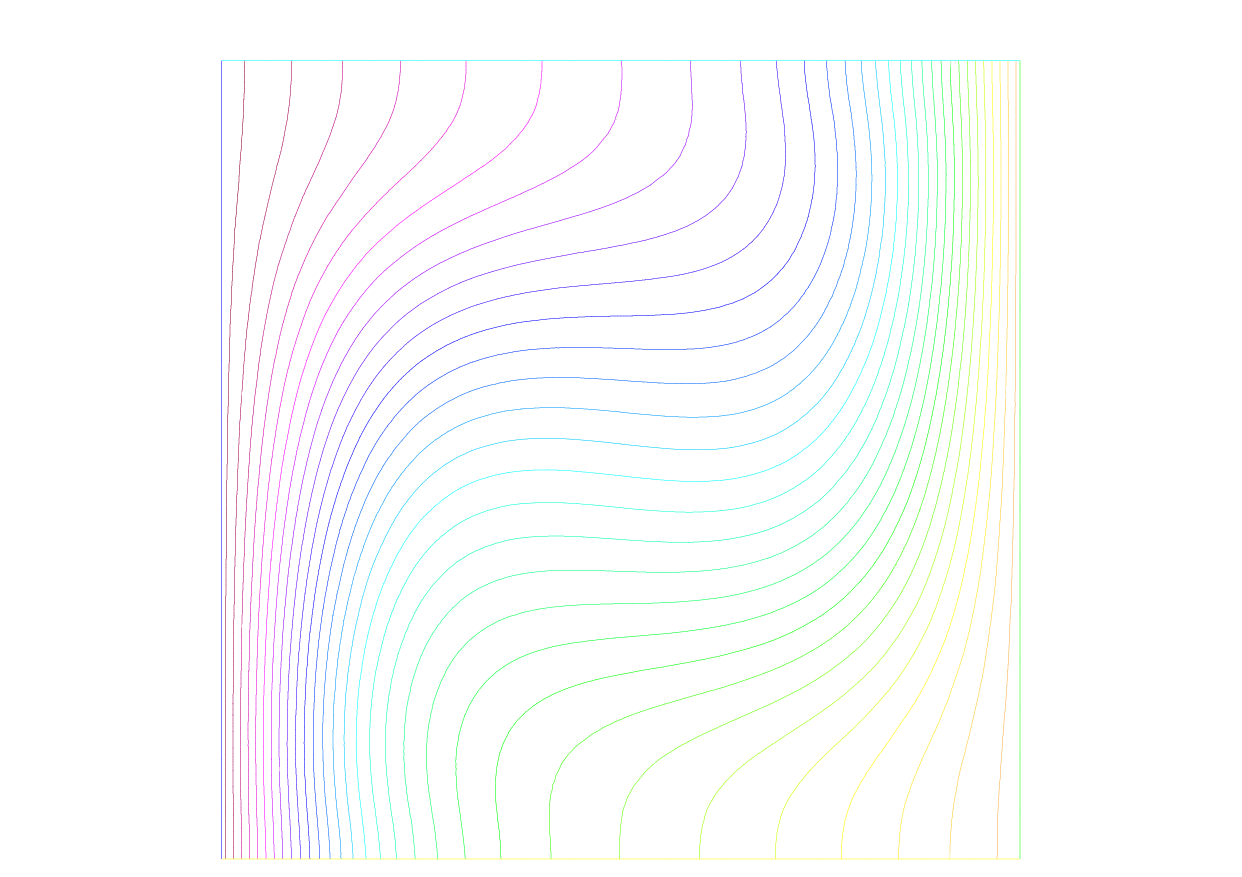}
\caption{Temperature }
\end{subfigure}
\begin{subfigure}{0.3\textwidth}
\includegraphics[width=0.9\linewidth]{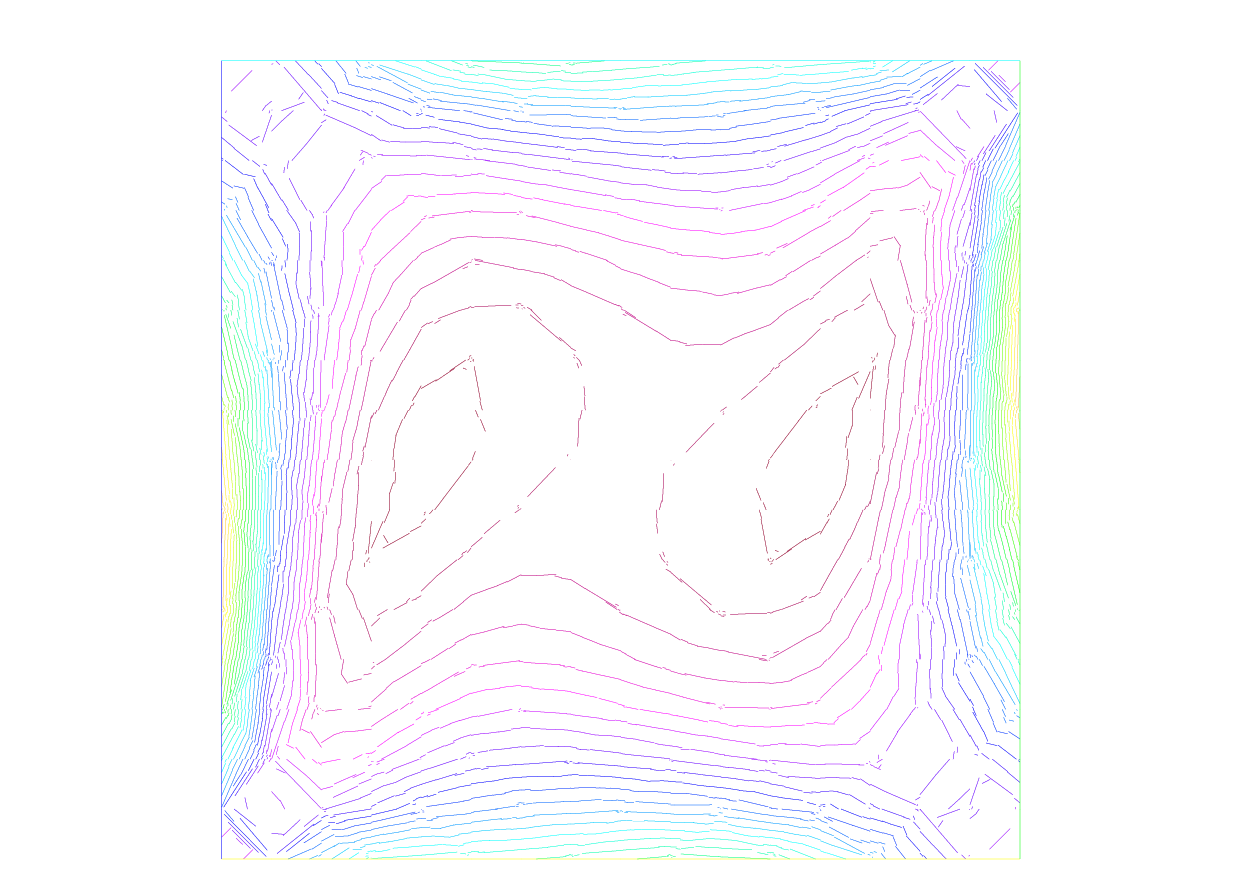}
\caption{Vorticity}
\end{subfigure}
\caption{$Ra=10^4$ Solution of natural convection without Coriolis.}
\label{fig:image1}
\end{figure}

\begin{figure}[H]
\begin{subfigure}{0.3\textwidth}
\includegraphics[width=0.9\linewidth]{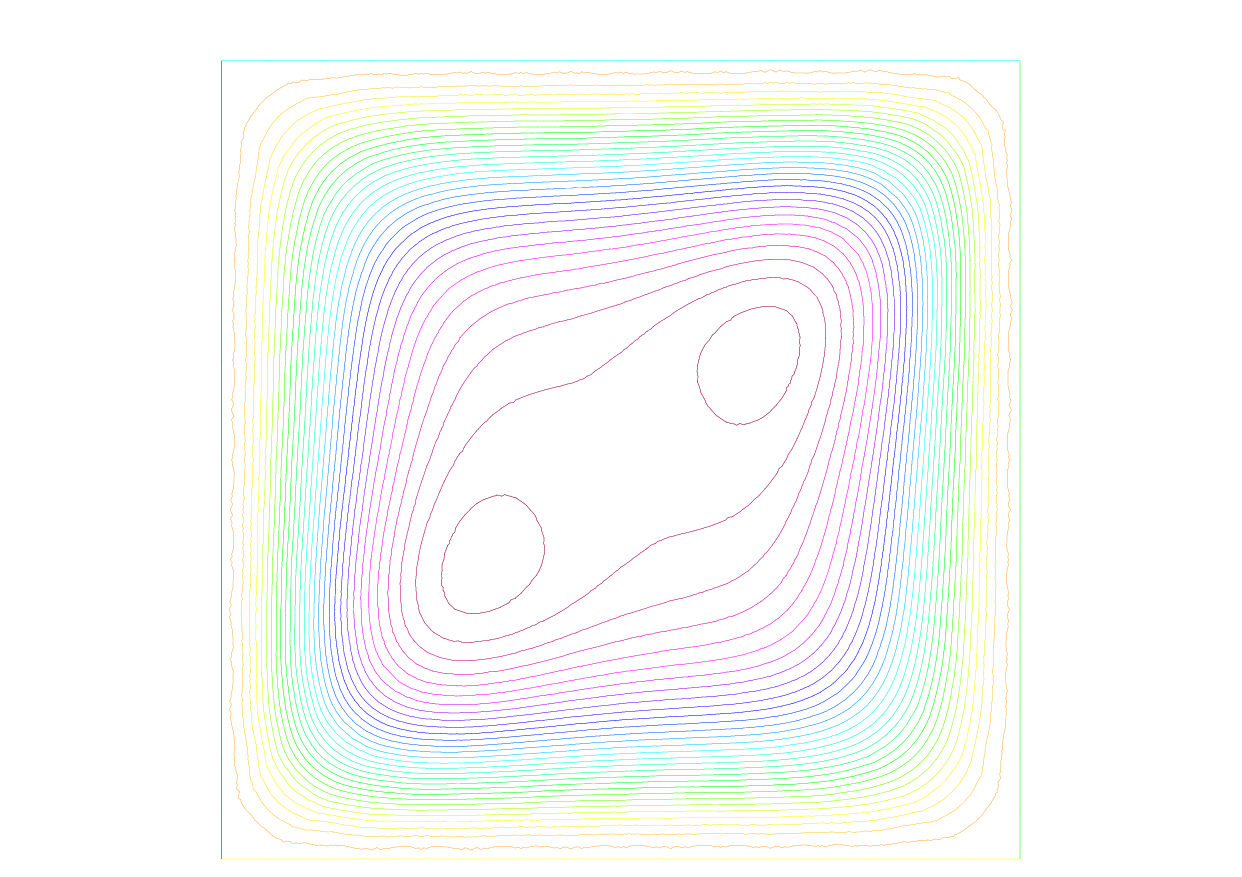} 
\caption{Streamlines}
\end{subfigure}
\begin{subfigure}{0.3\textwidth}
\includegraphics[width=0.9\linewidth]{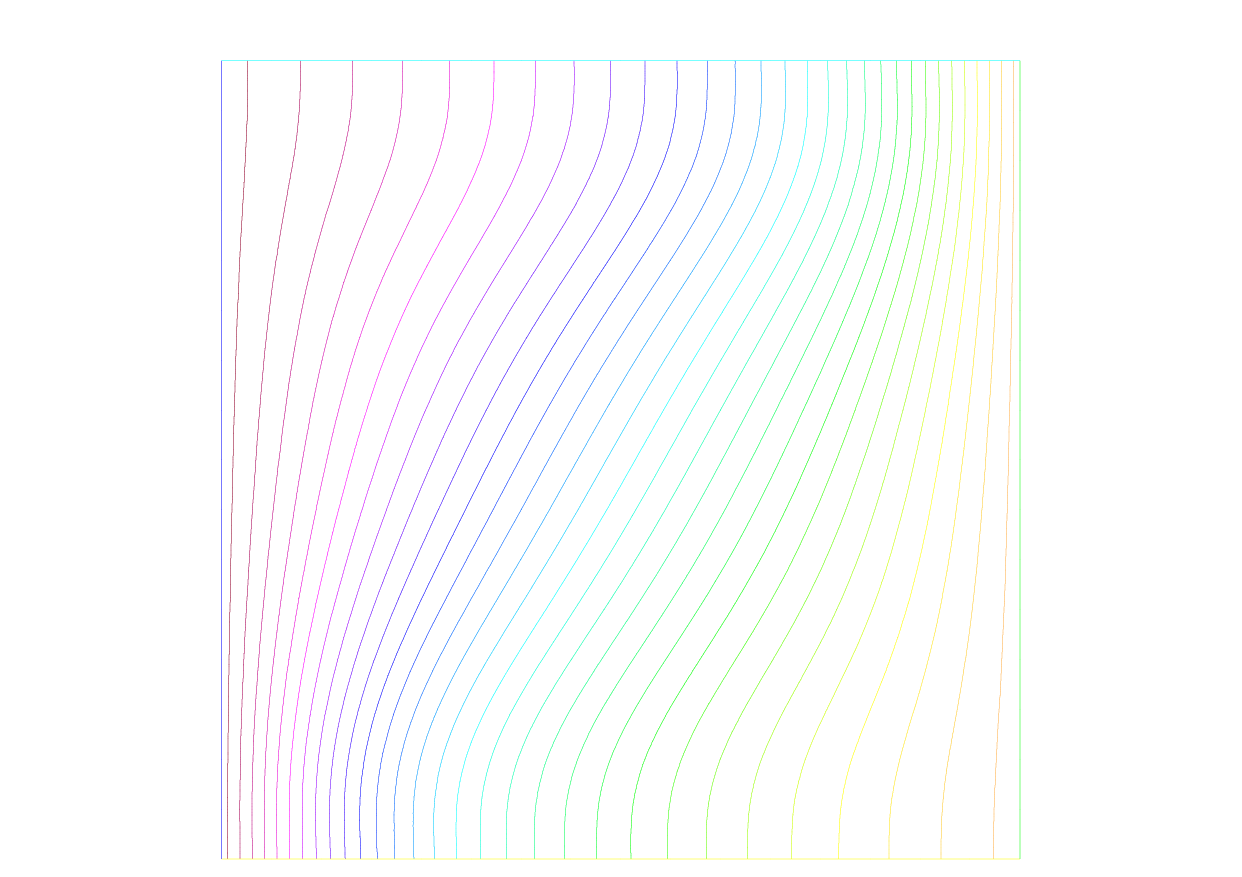}
\caption{Temperature }
\end{subfigure}
\begin{subfigure}{0.3\textwidth}
\includegraphics[width=0.9\linewidth]{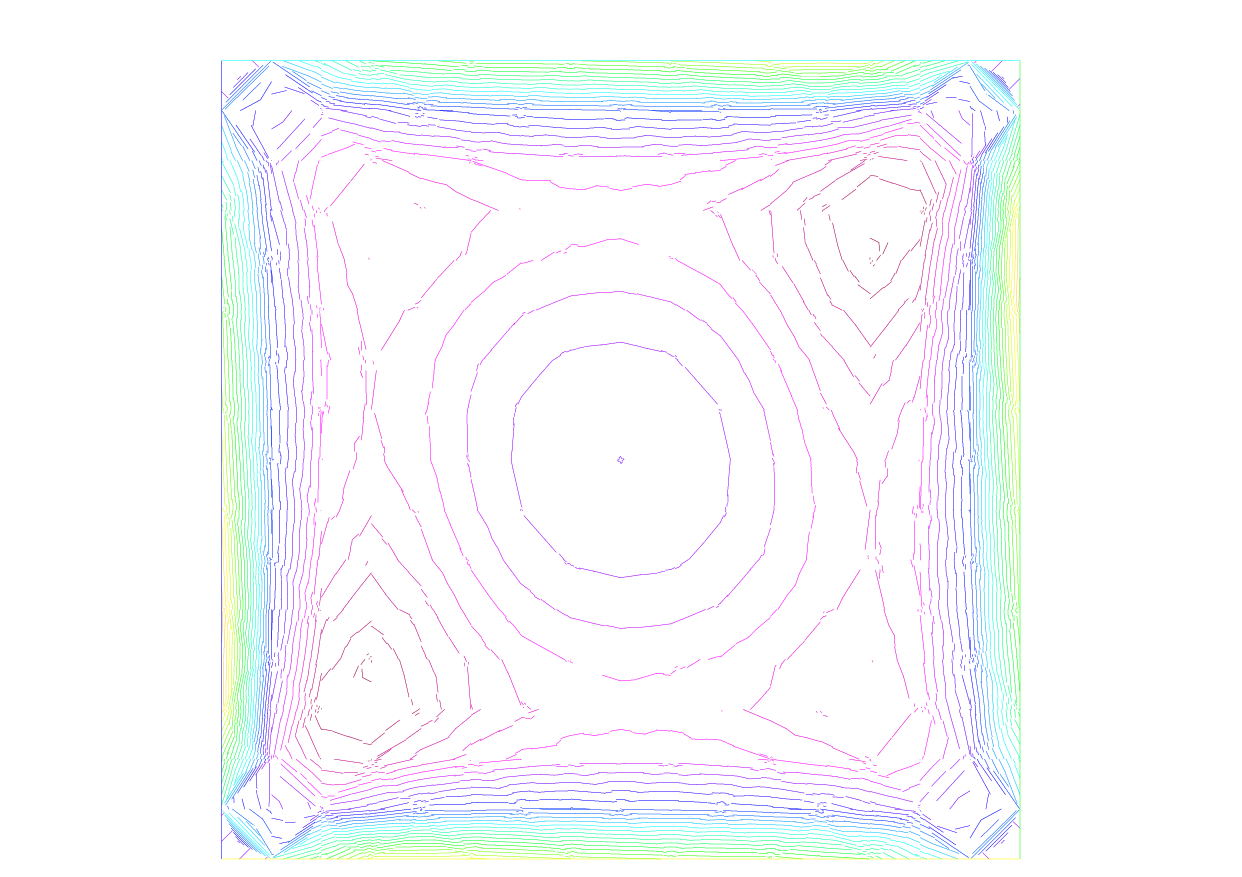}
\caption{Vorticity}
\end{subfigure}
\caption{$Ra=10^4$ Solution of natural convection with Coriolis.}
\label{fig:image2}
\end{figure}

We first depict the streamlines, temperature isolines, and vorticity with and without the Coriolis term, obtained by direct numerical simulation, to see the effect of the Coriolis force on the system. We pick $\omega= 5 \times 10^6$. Distortion caused by the Coriolis force is evident from Figure \ref{fig:image1} and Figure \ref{fig:image2}.

\begin{figure}[H]
\begin{subfigure}{0.3\textwidth}
\includegraphics[width=0.9\linewidth]{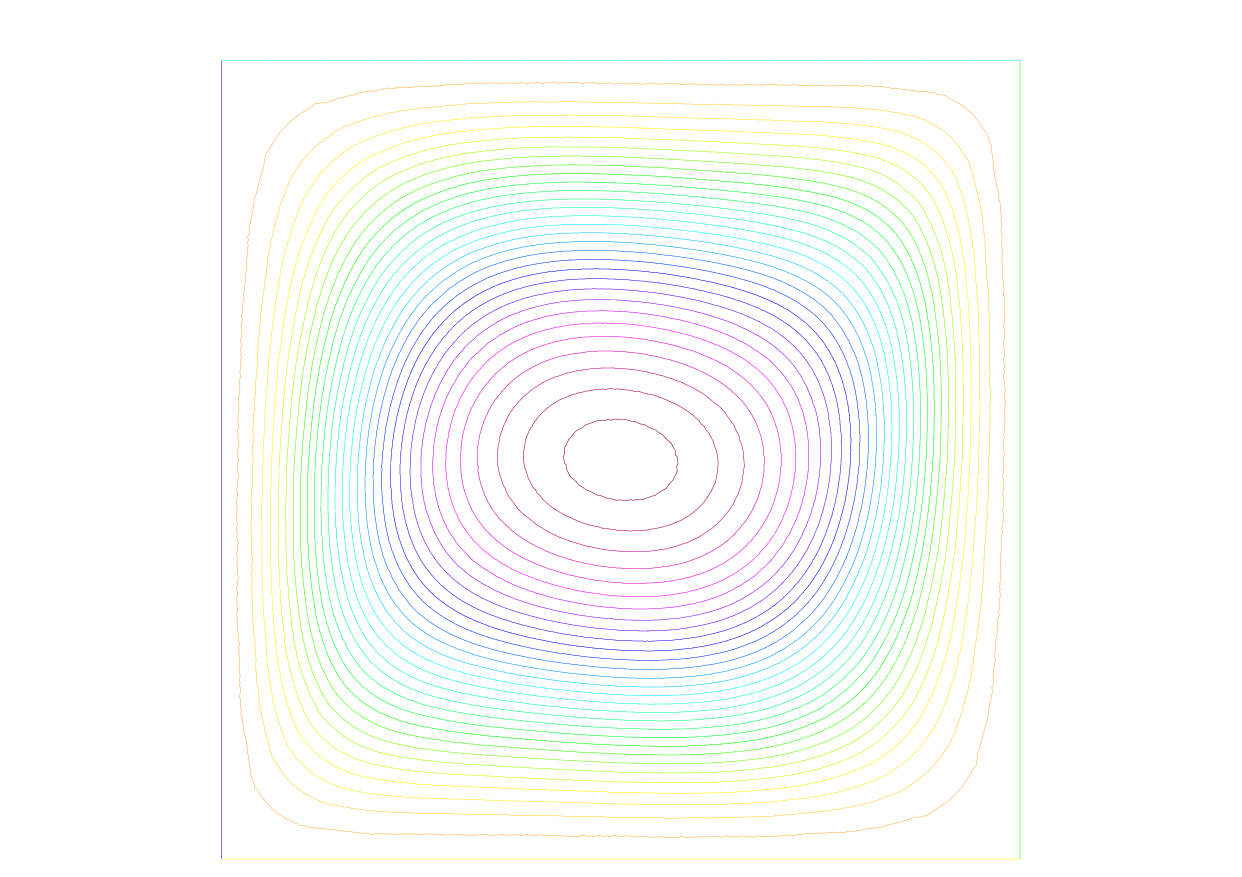} 
\caption{Streamlines}
\end{subfigure}
\begin{subfigure}{0.3\textwidth}
\includegraphics[width=0.9\linewidth]{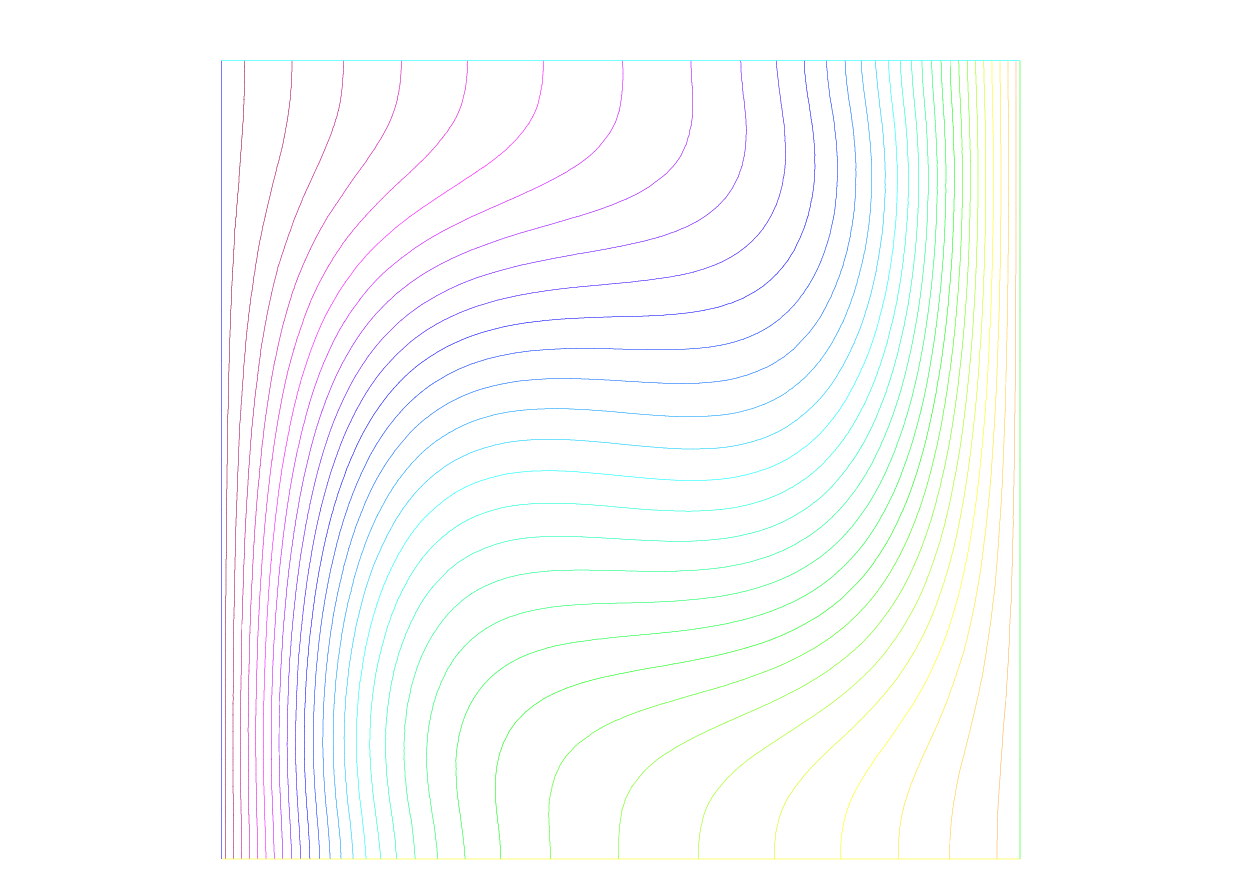}
\caption{Temperature }
\end{subfigure}
\begin{subfigure}{0.3\textwidth}
\includegraphics[width=0.9\linewidth]{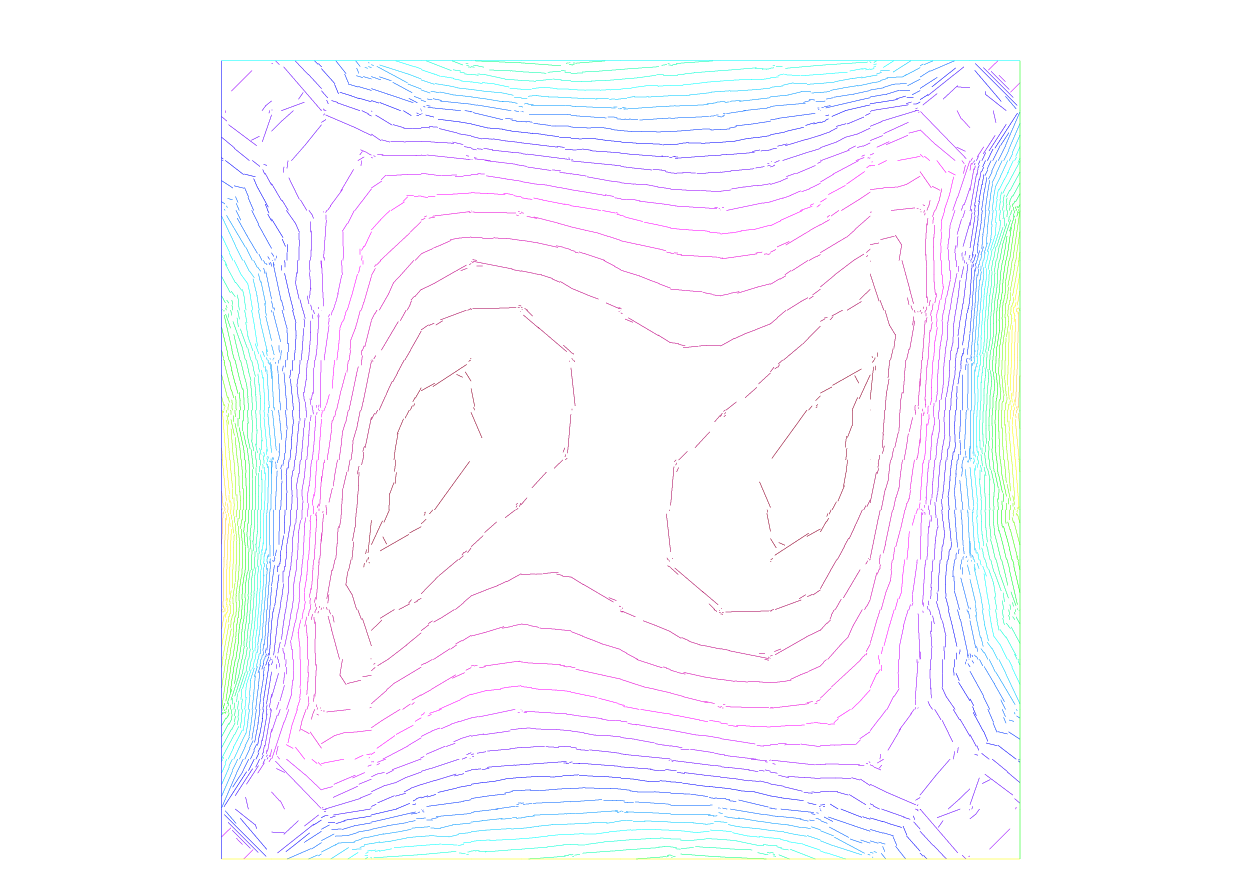}
\caption{Vorticity}
\end{subfigure}
\caption{$Ra=10^4$ Data assimilation solution for $\chi =1$.}
\label{fig:image3}
\end{figure}

\begin{figure}[H]
\begin{subfigure}{0.3\textwidth}
\includegraphics[width=0.9\linewidth]{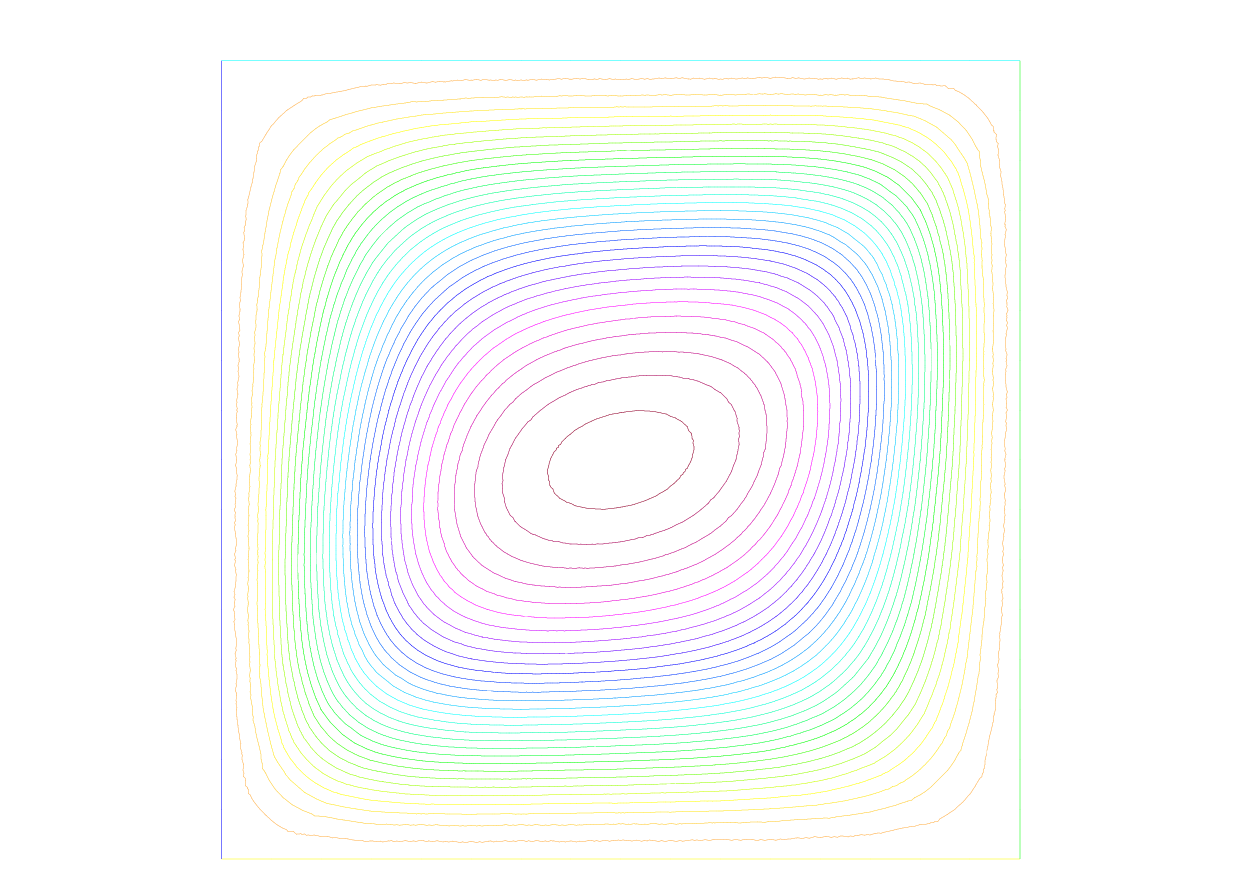} 
\caption{Streamlines}
\end{subfigure}
\begin{subfigure}{0.3\textwidth}
\includegraphics[width=0.9\linewidth]{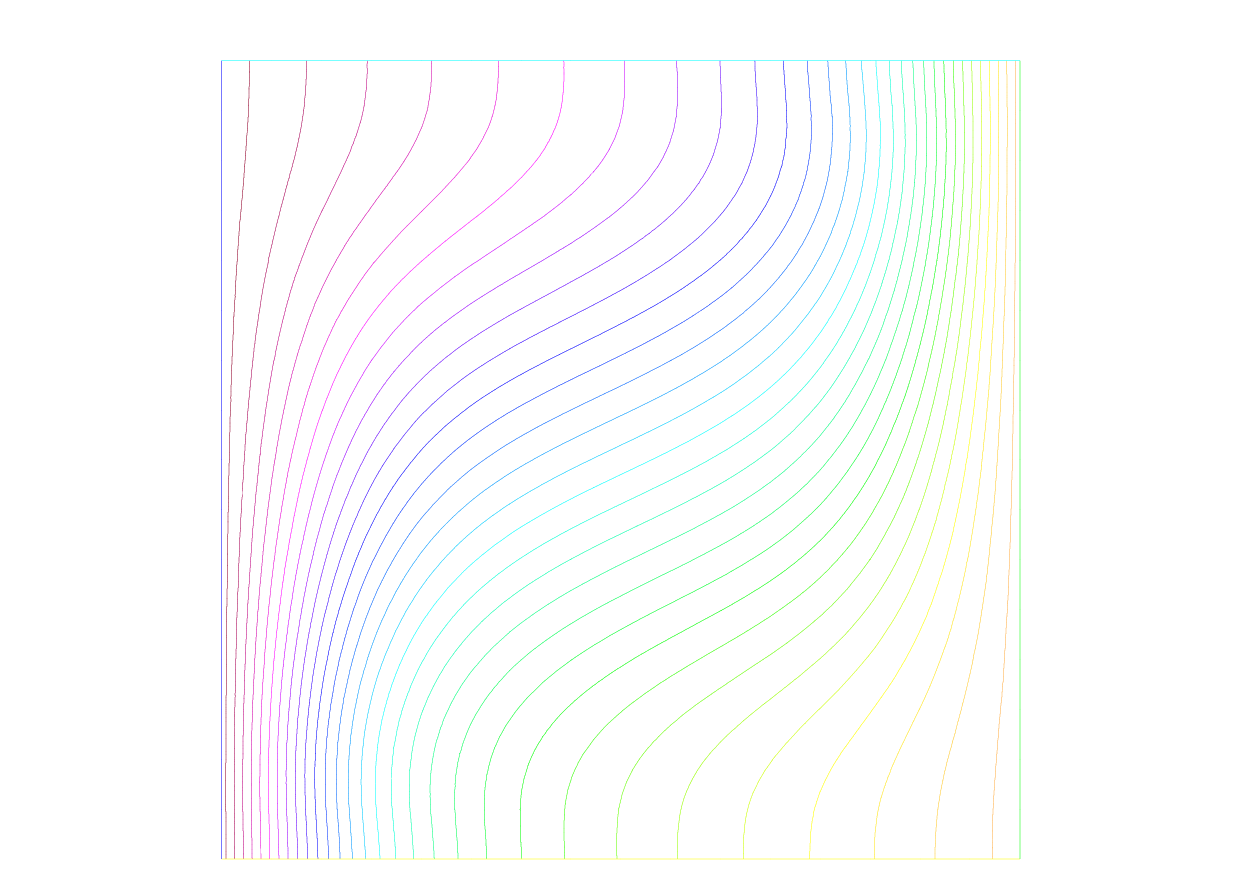}
\caption{Temperature }
\end{subfigure}
\begin{subfigure}{0.3\textwidth}
\includegraphics[width=0.9\linewidth]{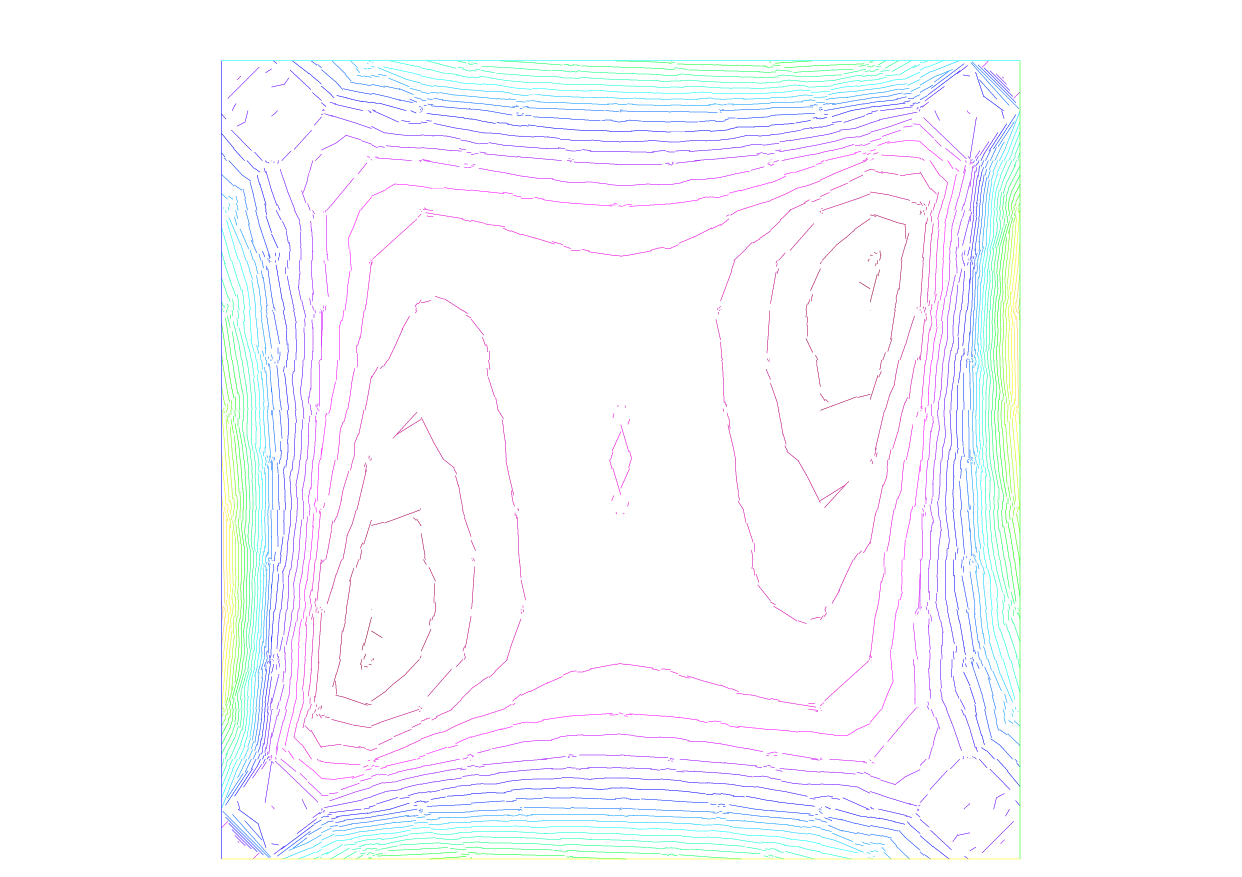}
\caption{Vorticity}
\end{subfigure}
\caption{$Ra=10^4$ Data assimilation solution for $\chi =10^2$.}
\label{fig:image4}
\end{figure}

\begin{figure}[H]
\begin{subfigure}{0.3\textwidth}
\includegraphics[width=0.9\linewidth]{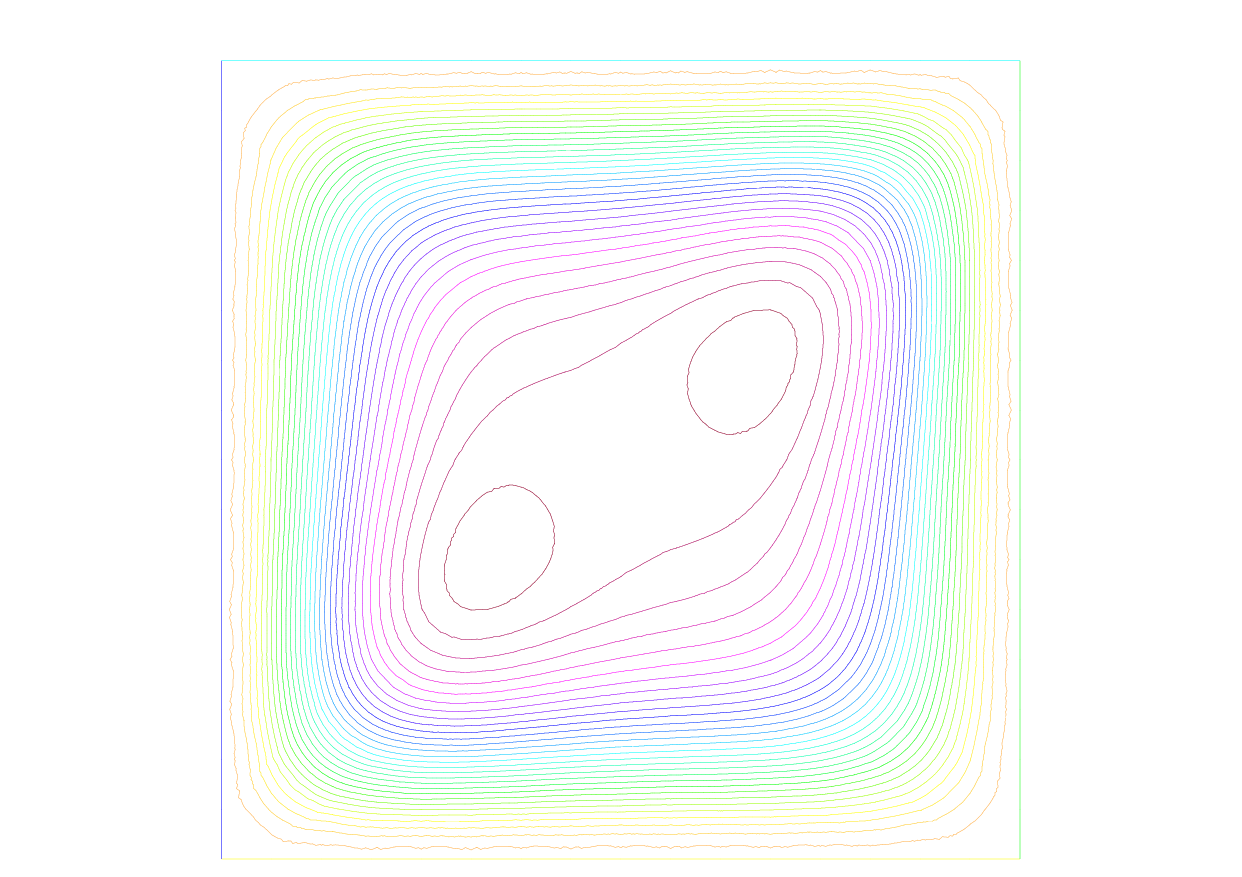} 
\caption{Streamlines}
\end{subfigure}
\begin{subfigure}{0.3\textwidth}
\includegraphics[width=0.9\linewidth]{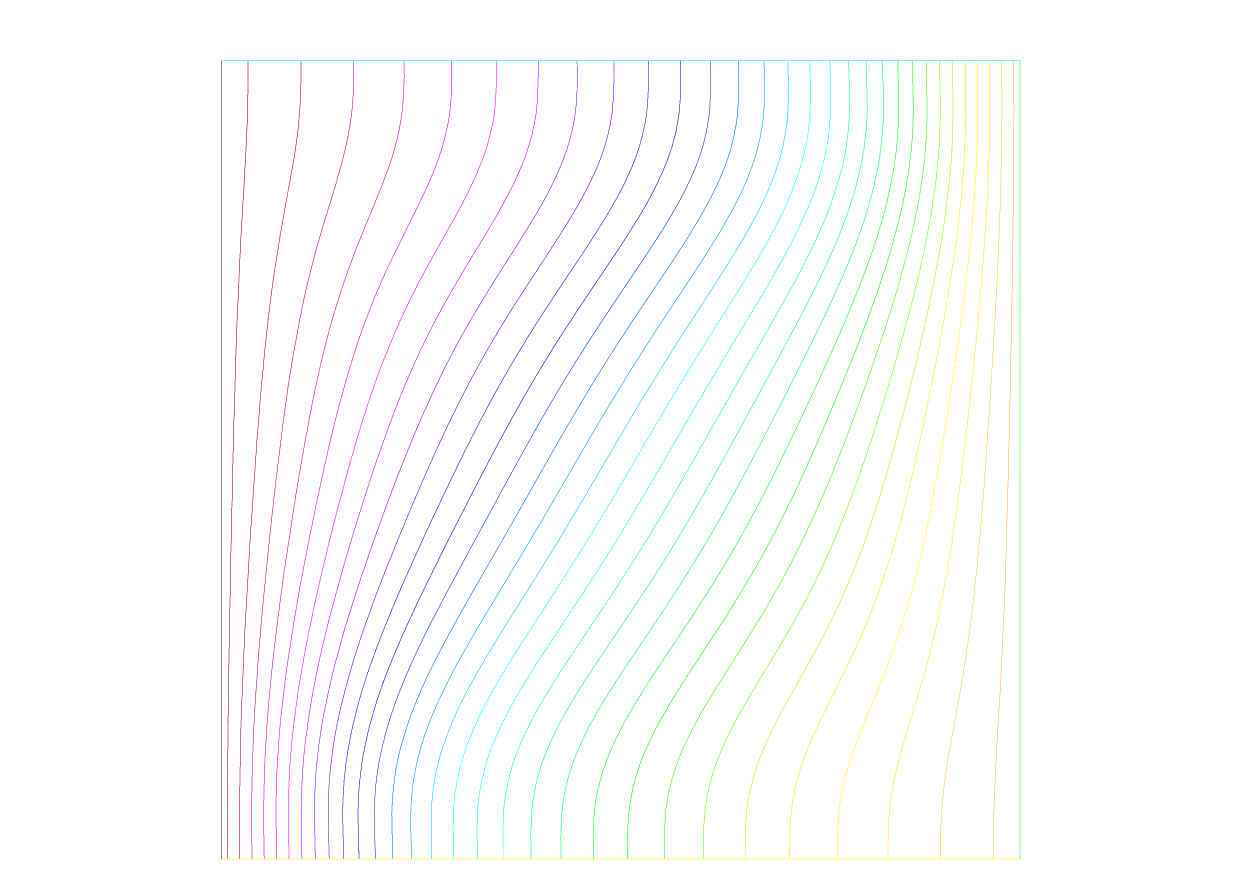}
\caption{Temperature }
\end{subfigure}
\begin{subfigure}{0.3\textwidth}
\includegraphics[width=0.9\linewidth]{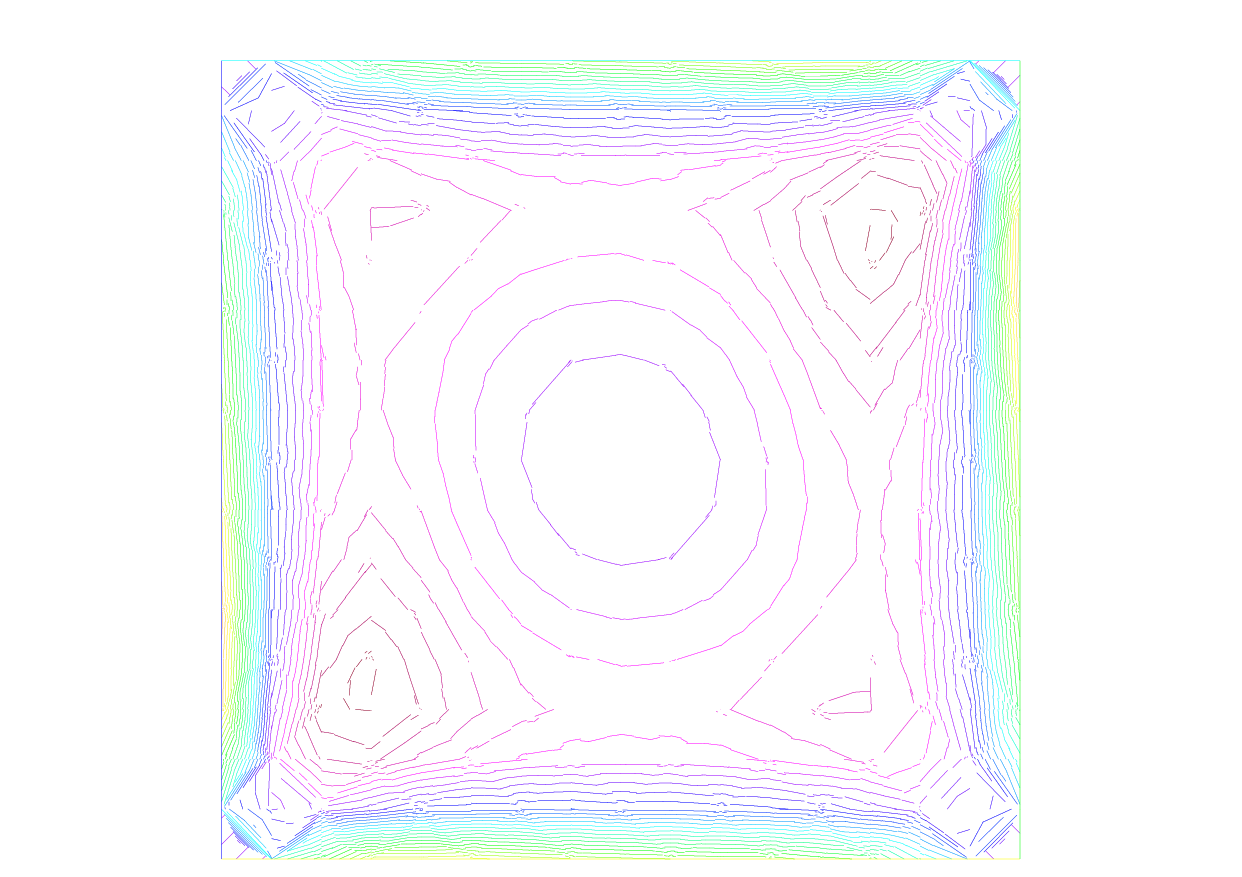}
\caption{Vorticity}
\end{subfigure}
\caption{$Ra=10^4$ Data assimilation solution for $\chi =10^4$.}
\label{fig:image5}
\end{figure}

\begin{figure}[H]
\begin{subfigure}{0.3\textwidth}
\includegraphics[width=0.9\linewidth]{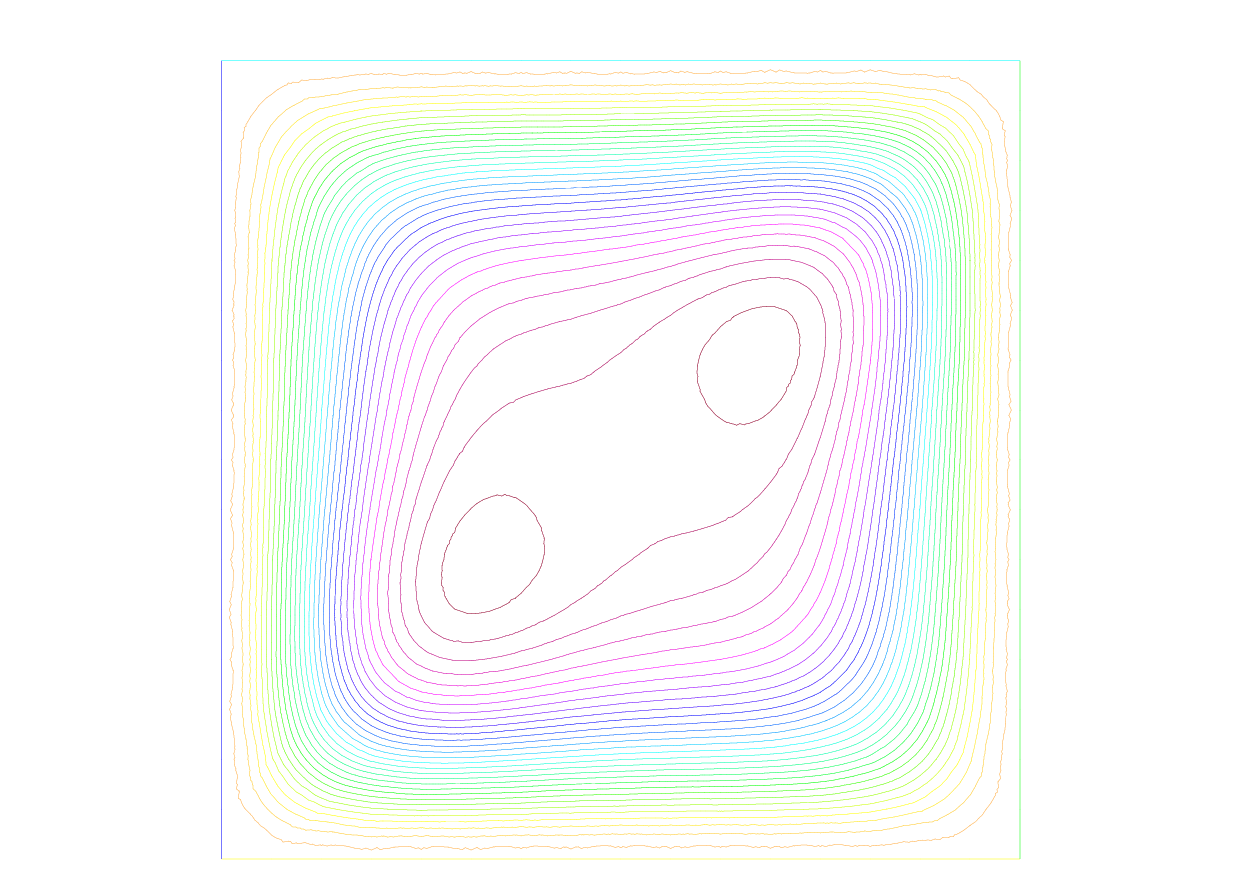} 
\caption{Streamlines}
\end{subfigure}
\begin{subfigure}{0.3\textwidth}
\includegraphics[width=0.9\linewidth]{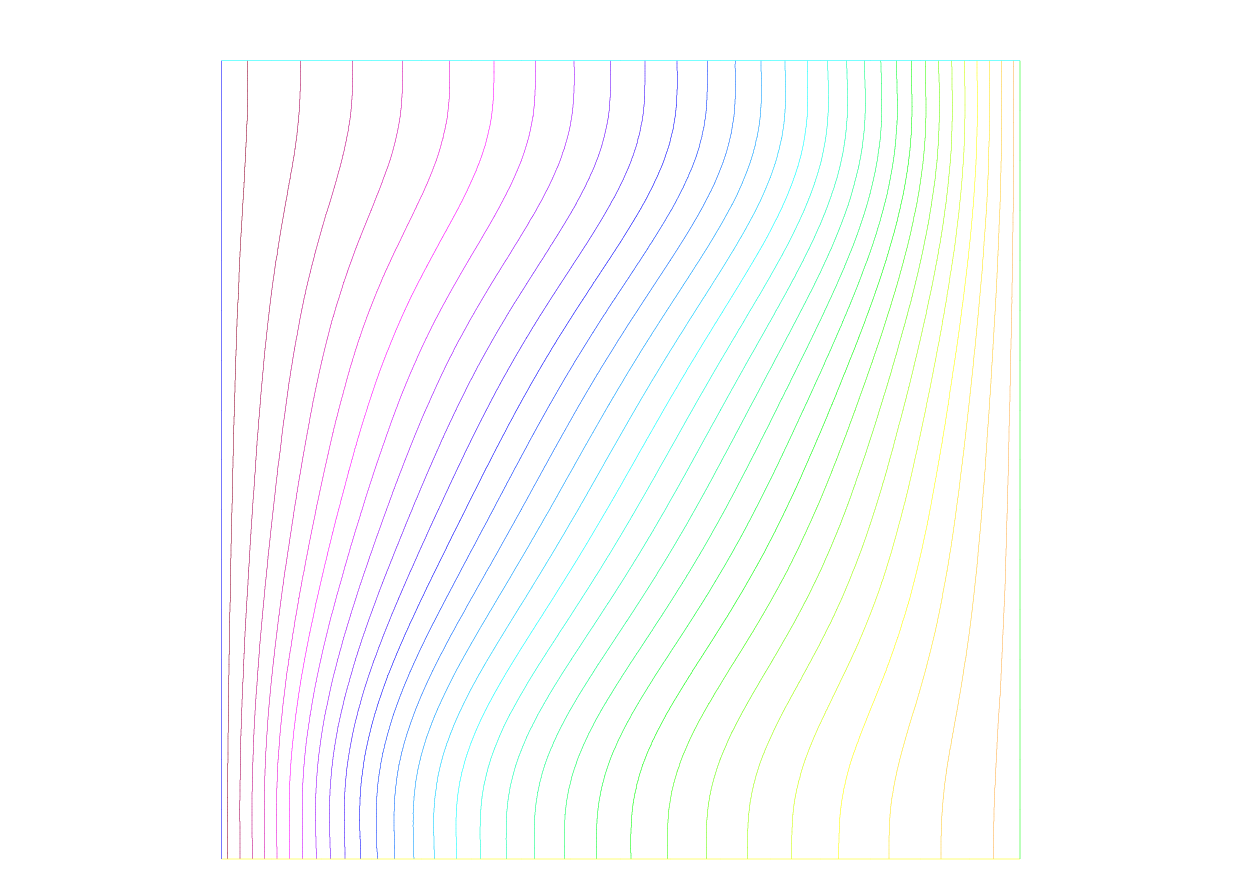}
\caption{Temperature }
\end{subfigure}
\begin{subfigure}{0.3\textwidth}
\includegraphics[width=0.9\linewidth]{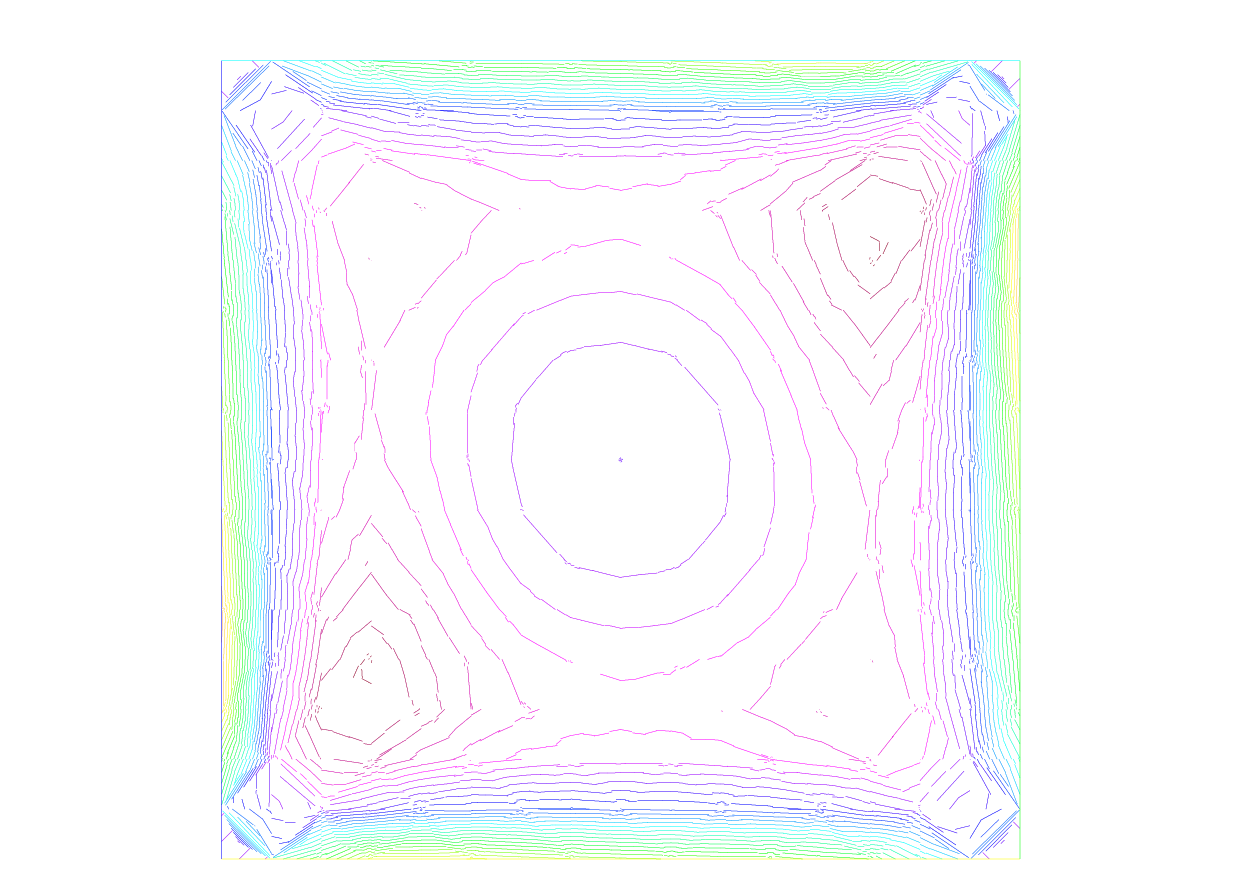}
\caption{Vorticity}
\end{subfigure}
\caption{$Ra=10^4$ Data assimilation solution for $\chi =10^6$.}
\label{fig:image6}
\end{figure}

Now, to see the effect of data assimilation by nudging, we present the flow patterns from (6.1), (6.2), (6.3) (without rotations) for different $\chi$ values. As $\chi$ increases all of the flow patterns become similar to the ones given in figure \ref{fig:image2}, evident in Figures \ref{fig:image3}, \ref{fig:image4}, \ref{fig:image5}, and \ref{fig:image6}.

\subsubsection{\textbf{Case II: $Ra=10^5$}}
\begin{figure}[H]
\begin{subfigure}{0.3\textwidth}
\includegraphics[width=0.9\linewidth]{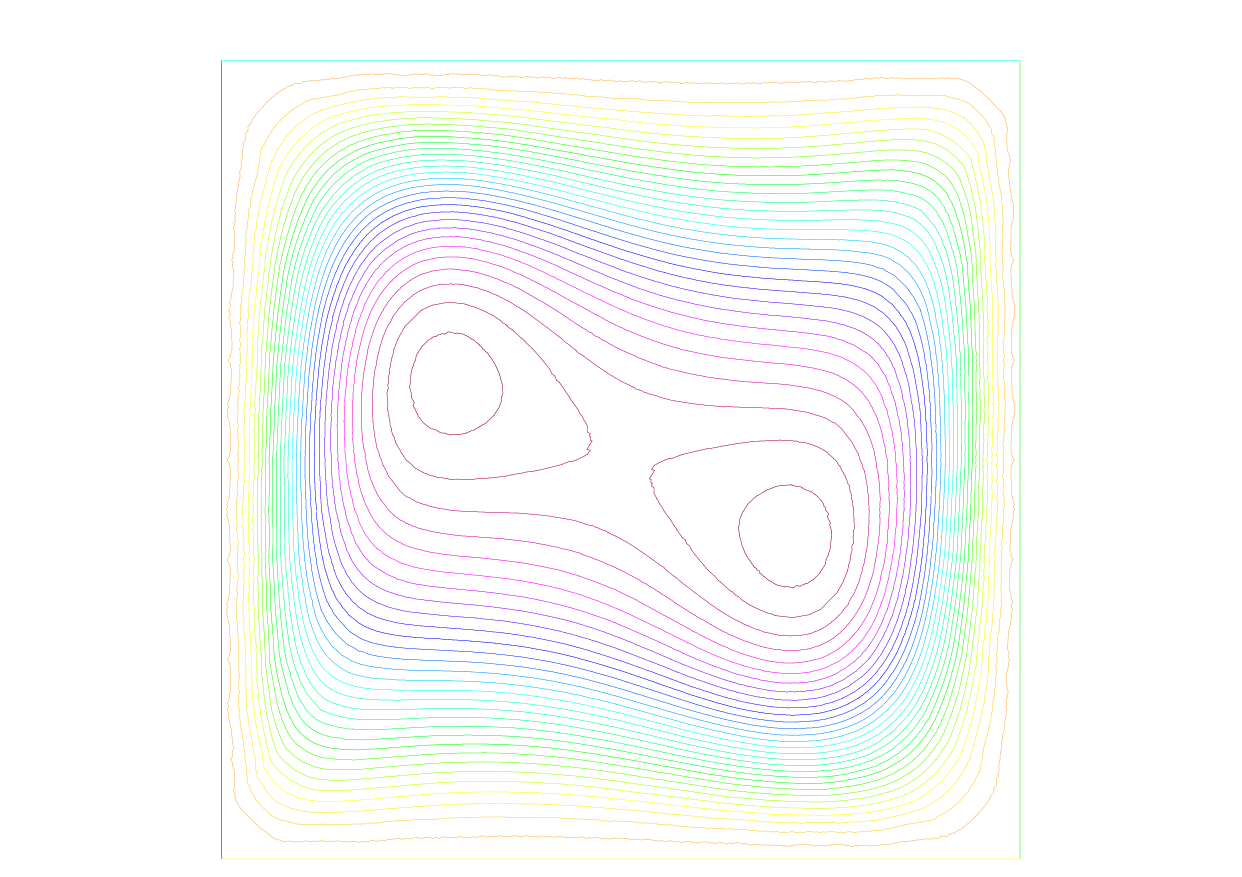} 
\caption{Streamlines}
\end{subfigure}
\begin{subfigure}{0.3\textwidth}
\includegraphics[width=0.9\linewidth]{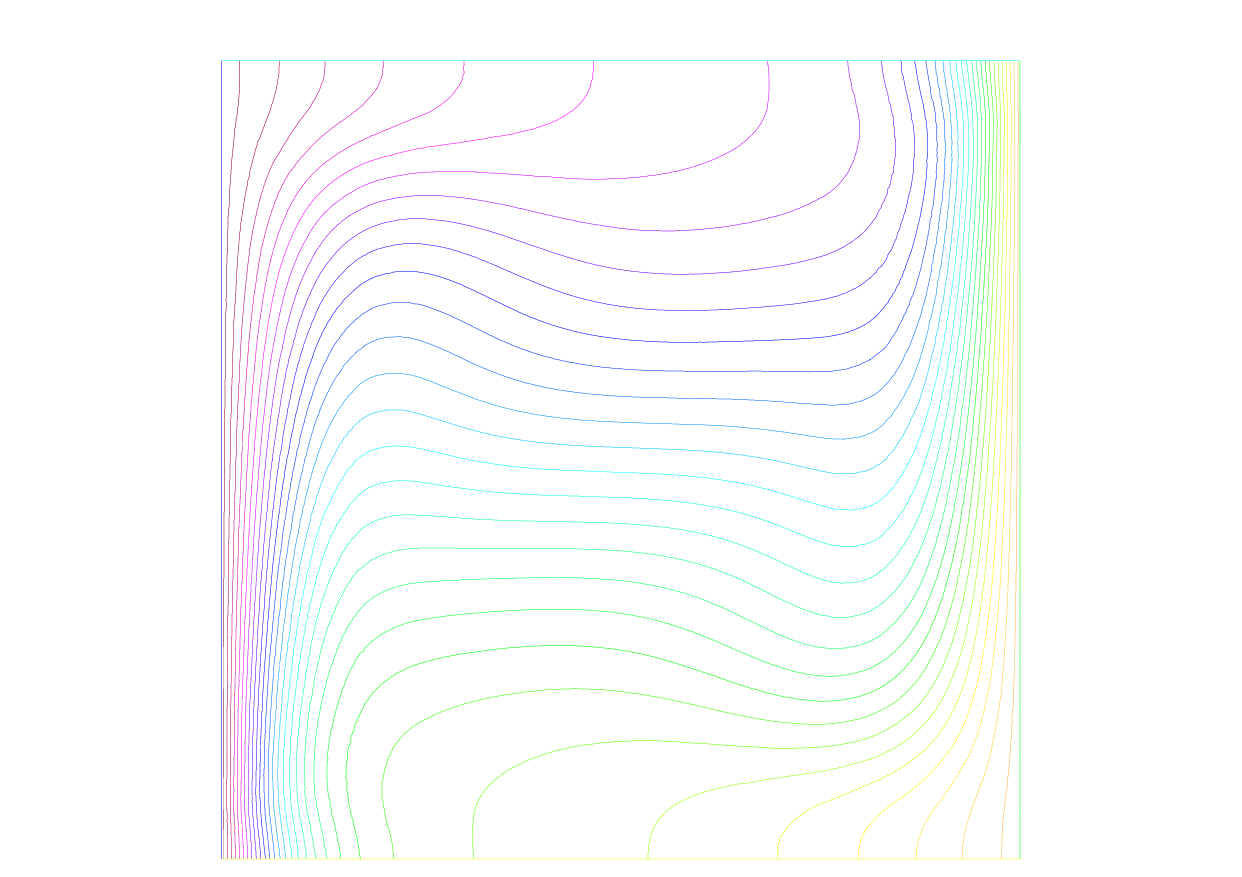}
\caption{Temperature }
\end{subfigure}
\begin{subfigure}{0.3\textwidth}
\includegraphics[width=0.9\linewidth]{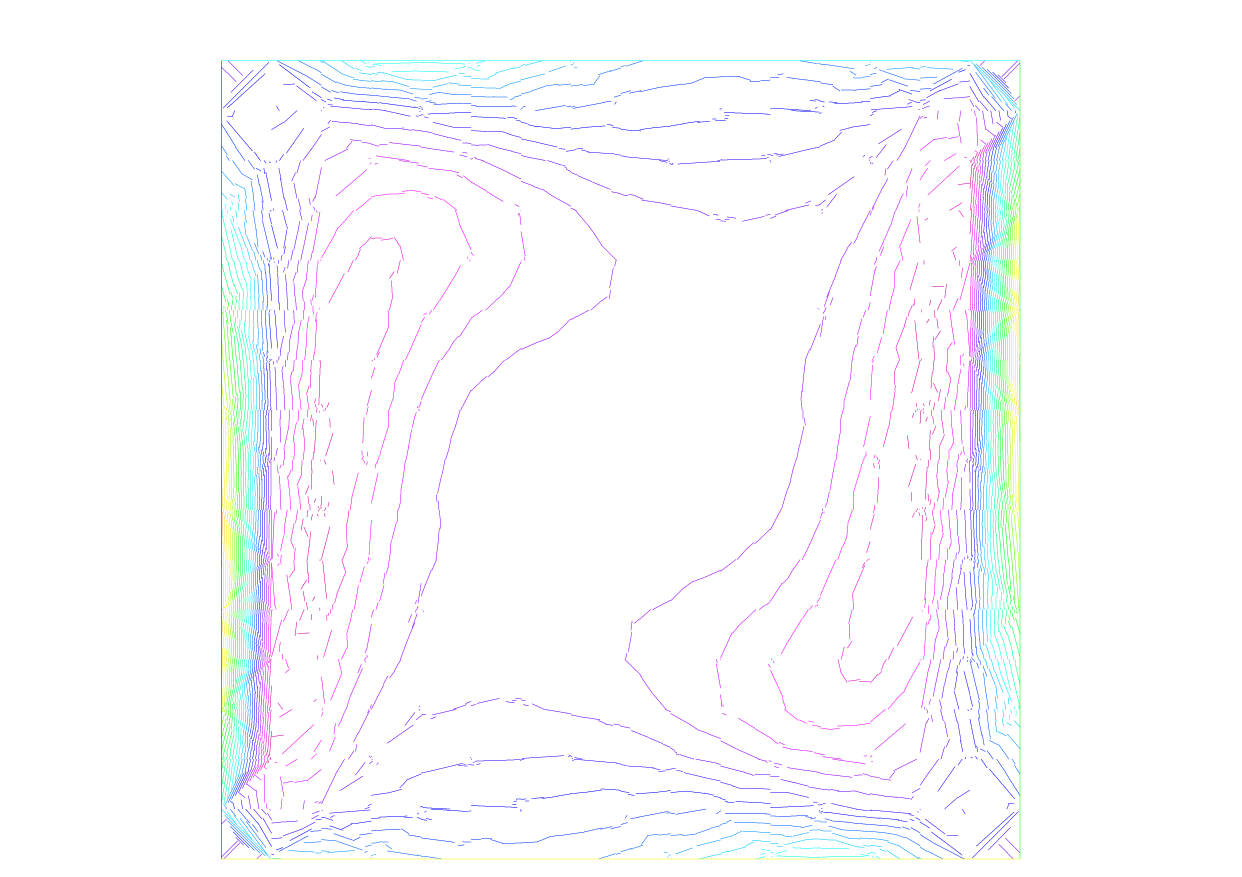}
\caption{Vorticity}
\end{subfigure}
\caption{$Ra=10^5$ Solution of Natural Convection without Coriolis.}
\label{fig:image7}
\end{figure}

\begin{figure}[H]
\begin{subfigure}{0.3\textwidth}
\includegraphics[width=0.9\linewidth]{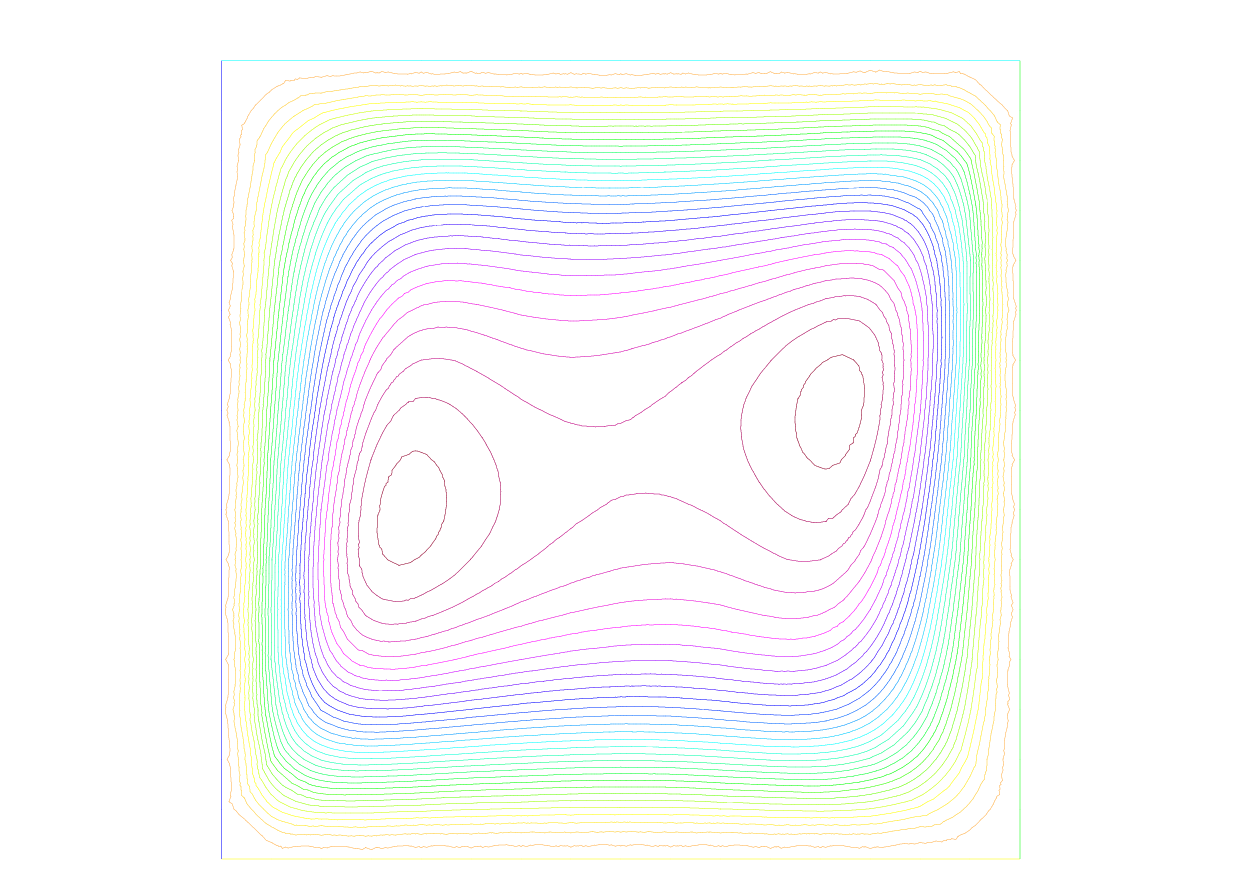} 
\caption{Streamlines}
\end{subfigure}
\begin{subfigure}{0.3\textwidth}
\includegraphics[width=0.9\linewidth]{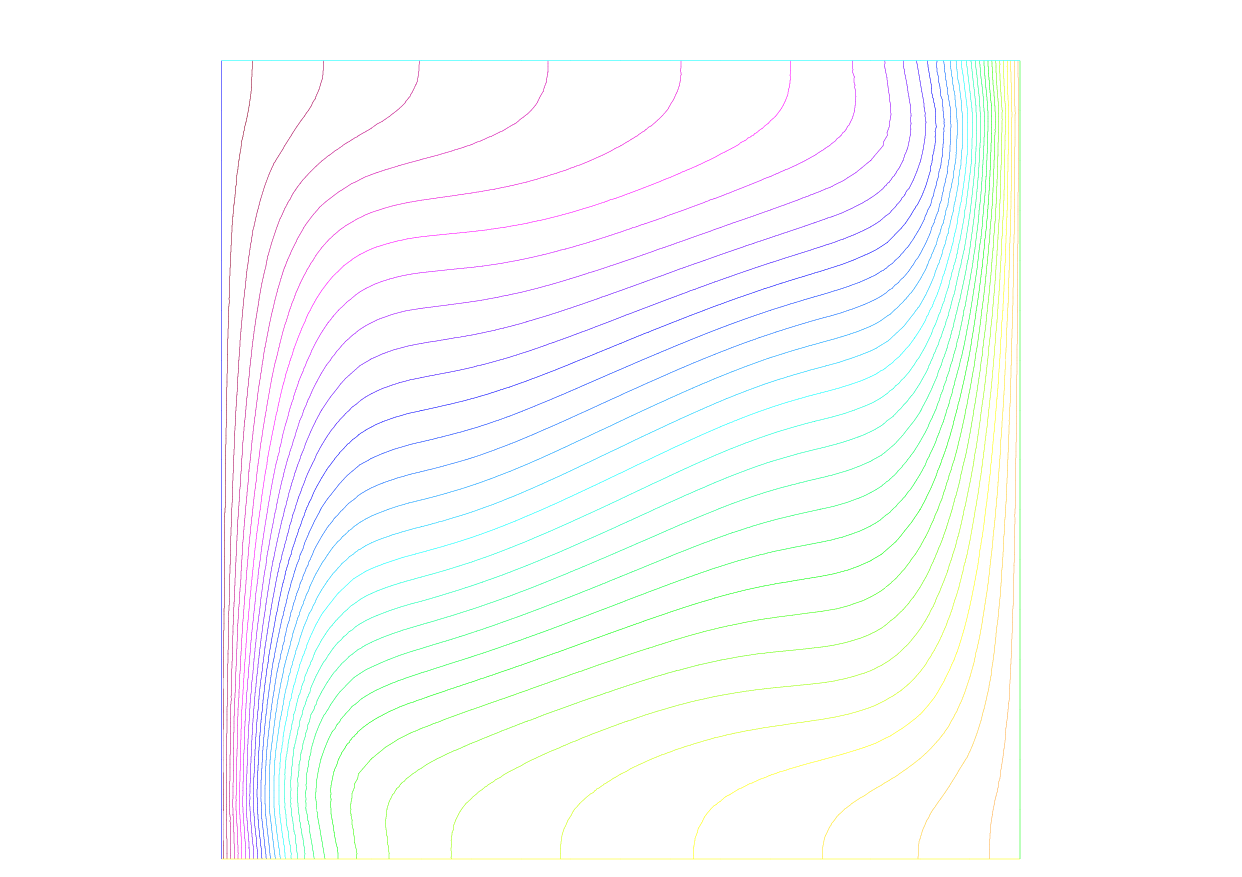}
\caption{Temperature }
\end{subfigure}
\begin{subfigure}{0.3\textwidth}
\includegraphics[width=0.9\linewidth]{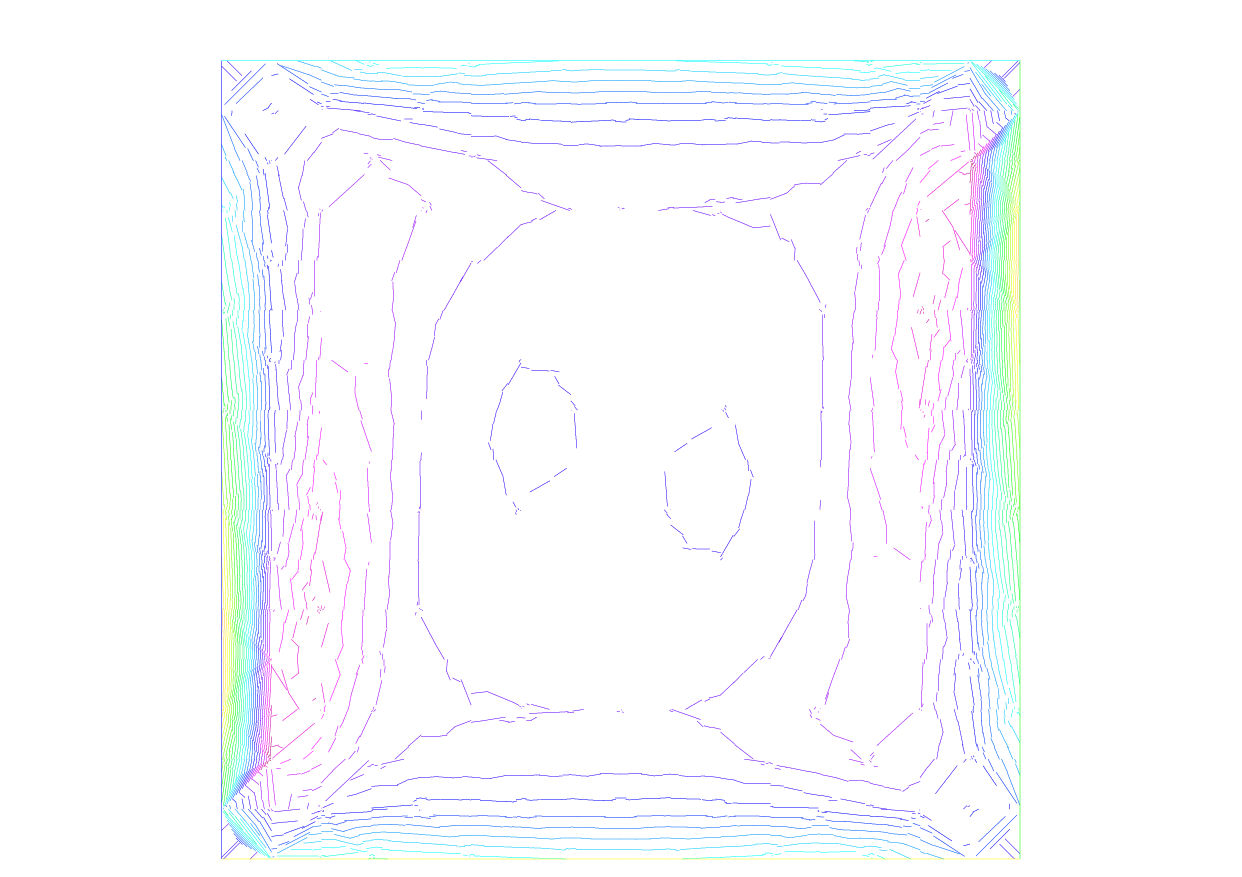}
\caption{Vorticity}
\end{subfigure}
\caption{$Ra=10^5$ Solution of Natural Convection with Coriolis.}
\label{fig:image8}
\end{figure}

We apply the same procedure as for the $Ra=10^4$ case. Solution graphics with and without Coriolis force are given in Figures \ref{fig:image7} and \ref{fig:image8}.

Now, we apply nudging and solve the data assimilation problem with various $\chi$ values.

\begin{figure}[H]
\begin{subfigure}{0.3\textwidth}
\includegraphics[width=0.9\linewidth]{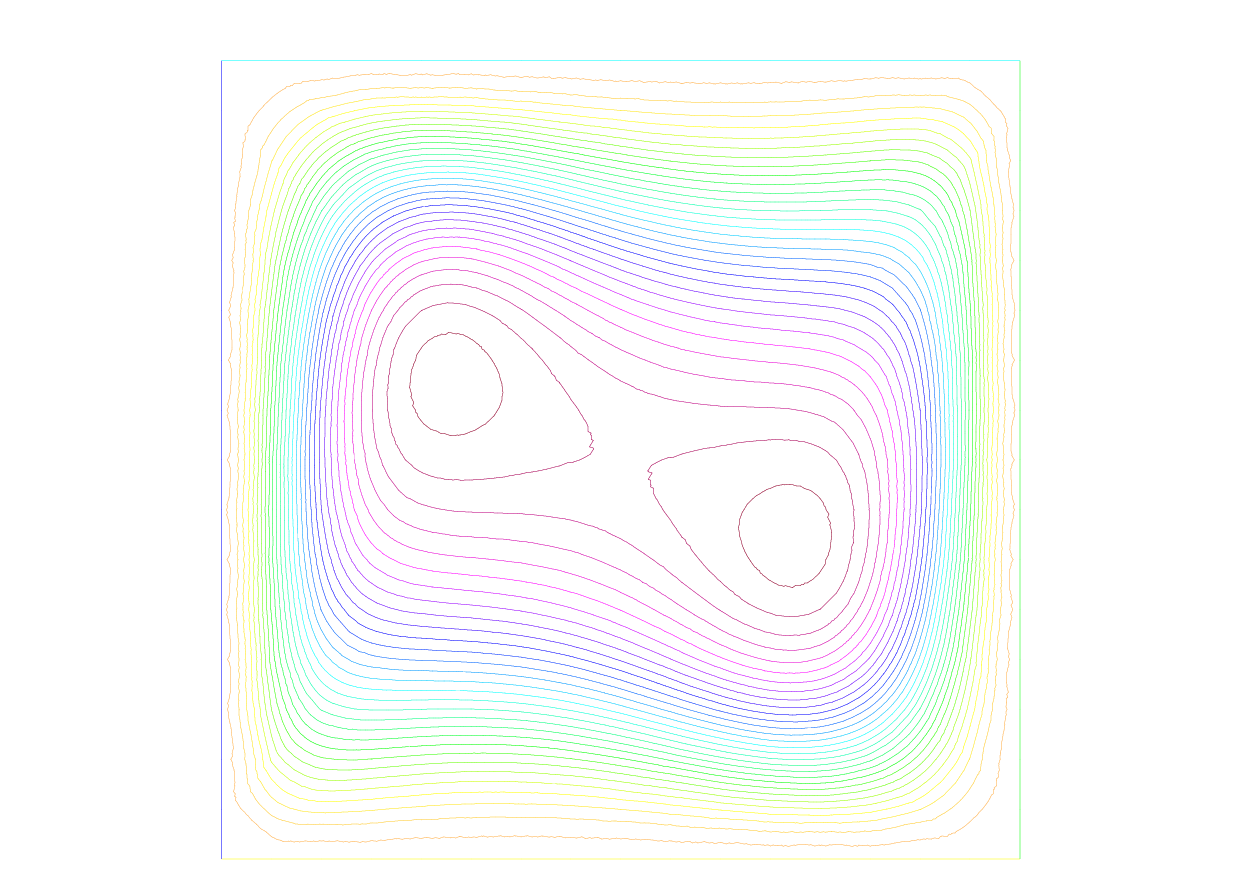} 
\caption{Streamlines}
\end{subfigure}
\begin{subfigure}{0.3\textwidth}
\includegraphics[width=0.9\linewidth]{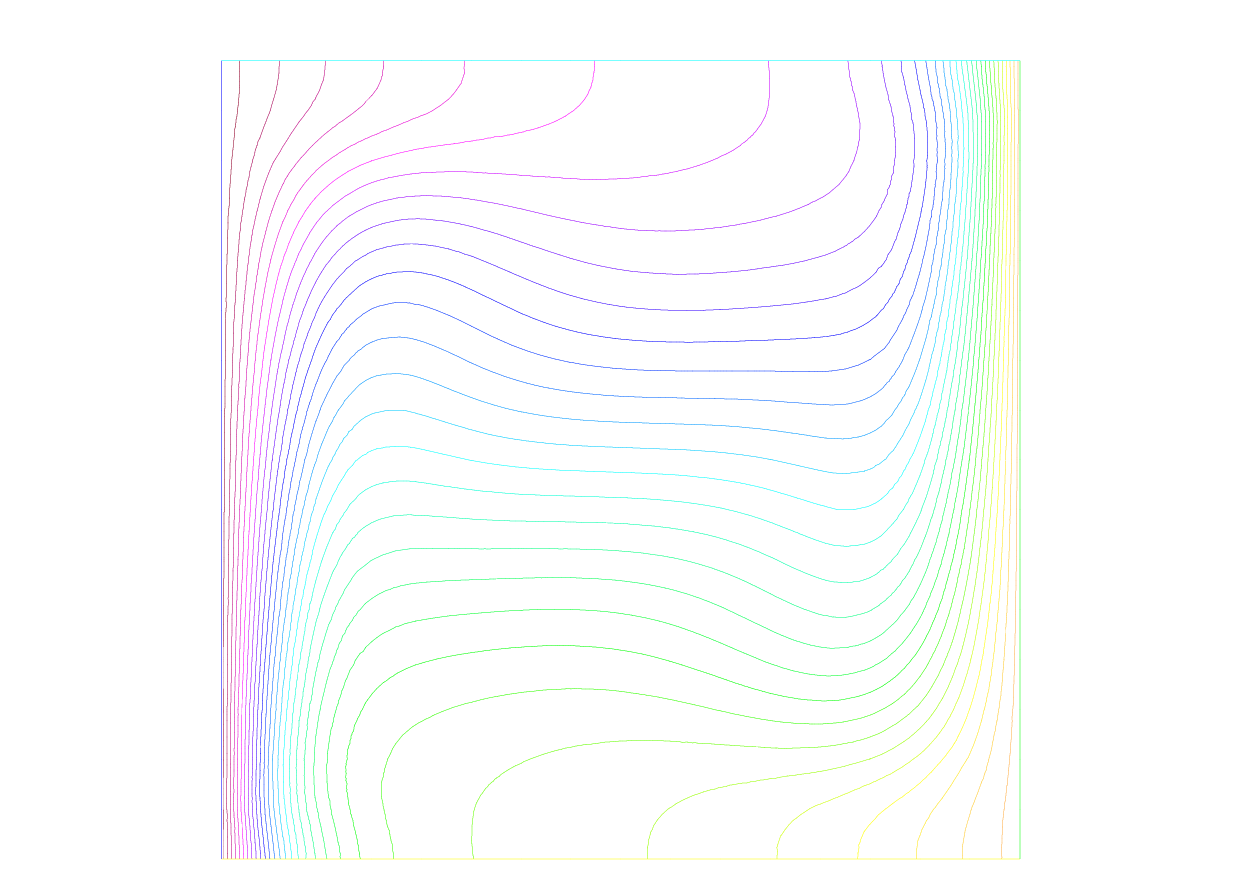}
\caption{Temperature }
\end{subfigure}
\begin{subfigure}{0.3\textwidth}
\includegraphics[width=0.9\linewidth]{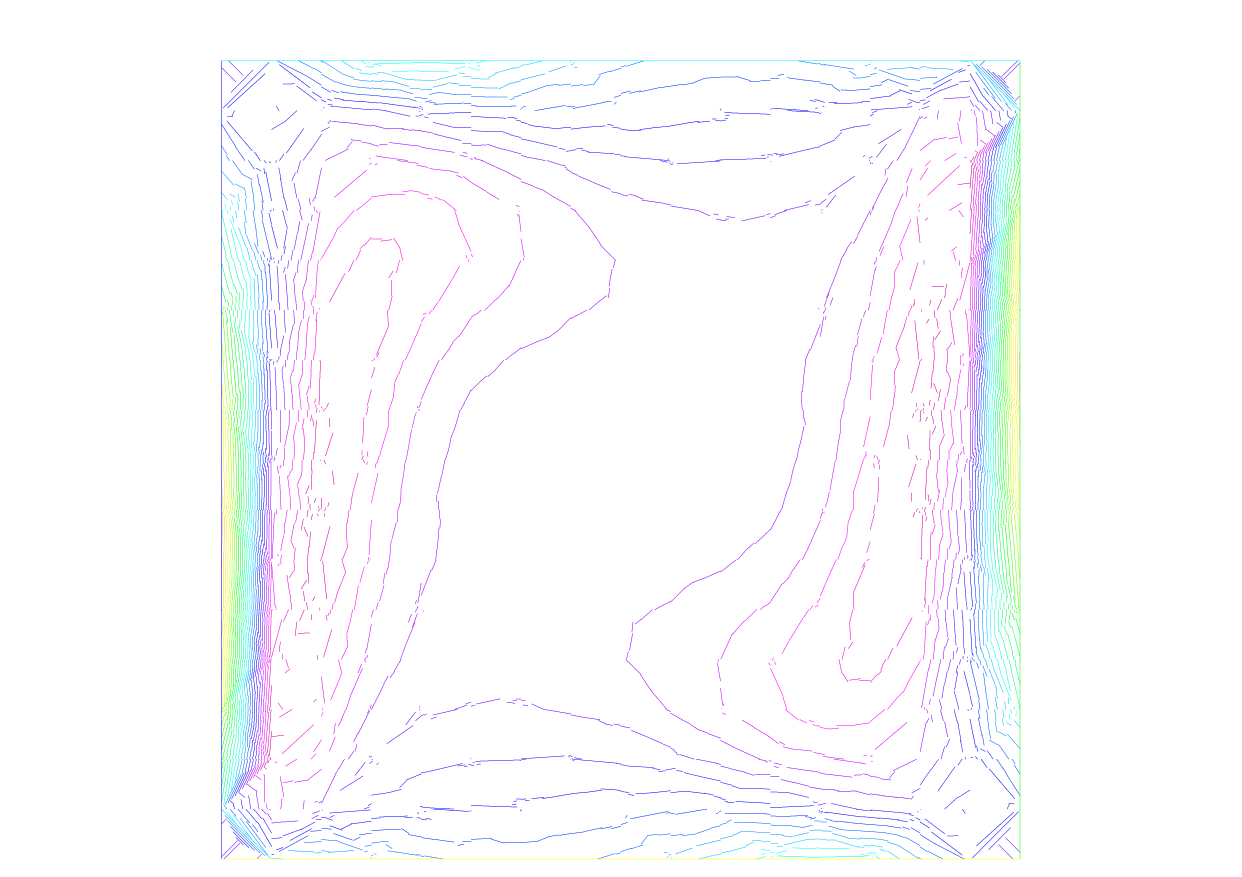}
\caption{Vorticity}
\end{subfigure}
\caption{$Ra=10^5$ Data Assimilation solution for $\chi =1$.}
\label{fig:image9}
\end{figure}

\begin{figure}[H]
\begin{subfigure}{0.3\textwidth}
\includegraphics[width=0.9\linewidth]{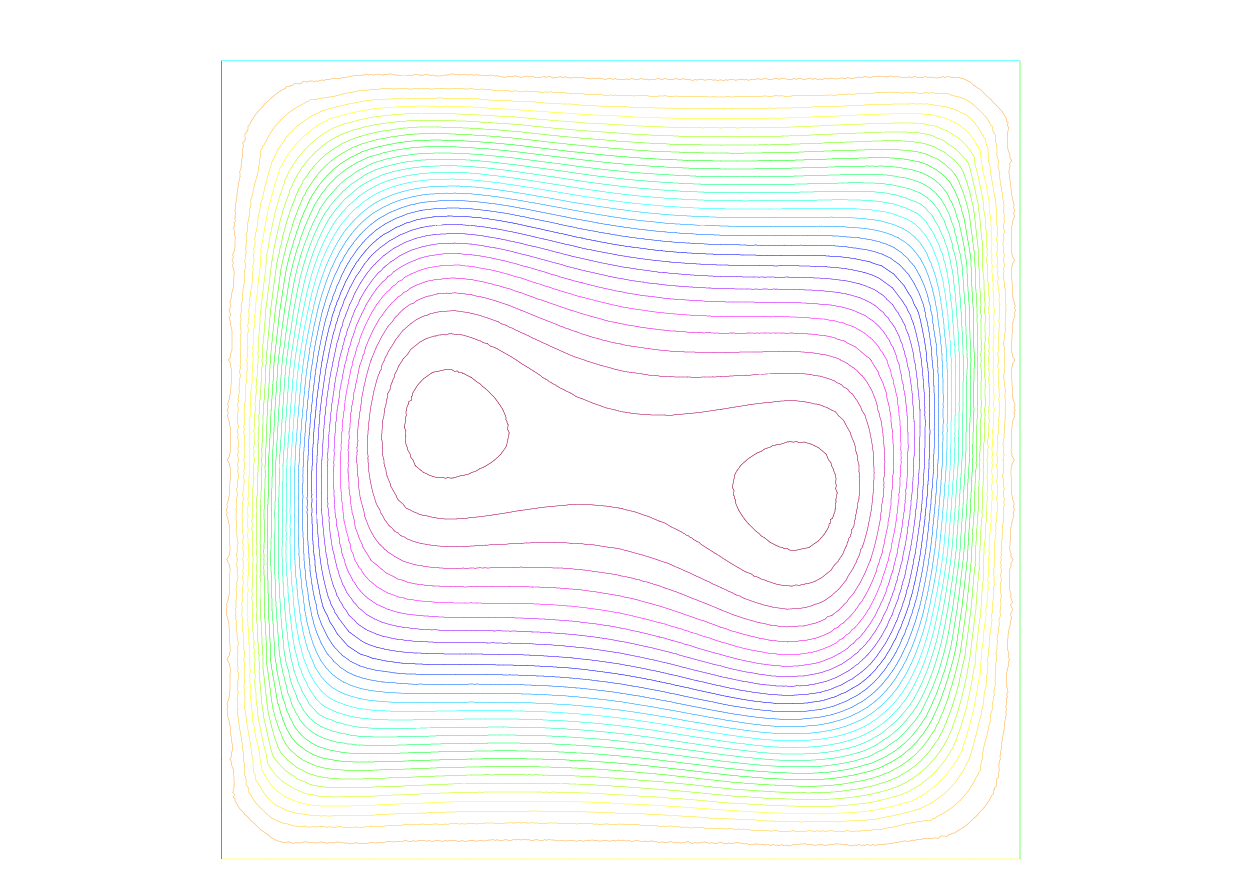} 
\caption{Streamlines}
\end{subfigure}
\begin{subfigure}{0.3\textwidth}
\includegraphics[width=0.9\linewidth]{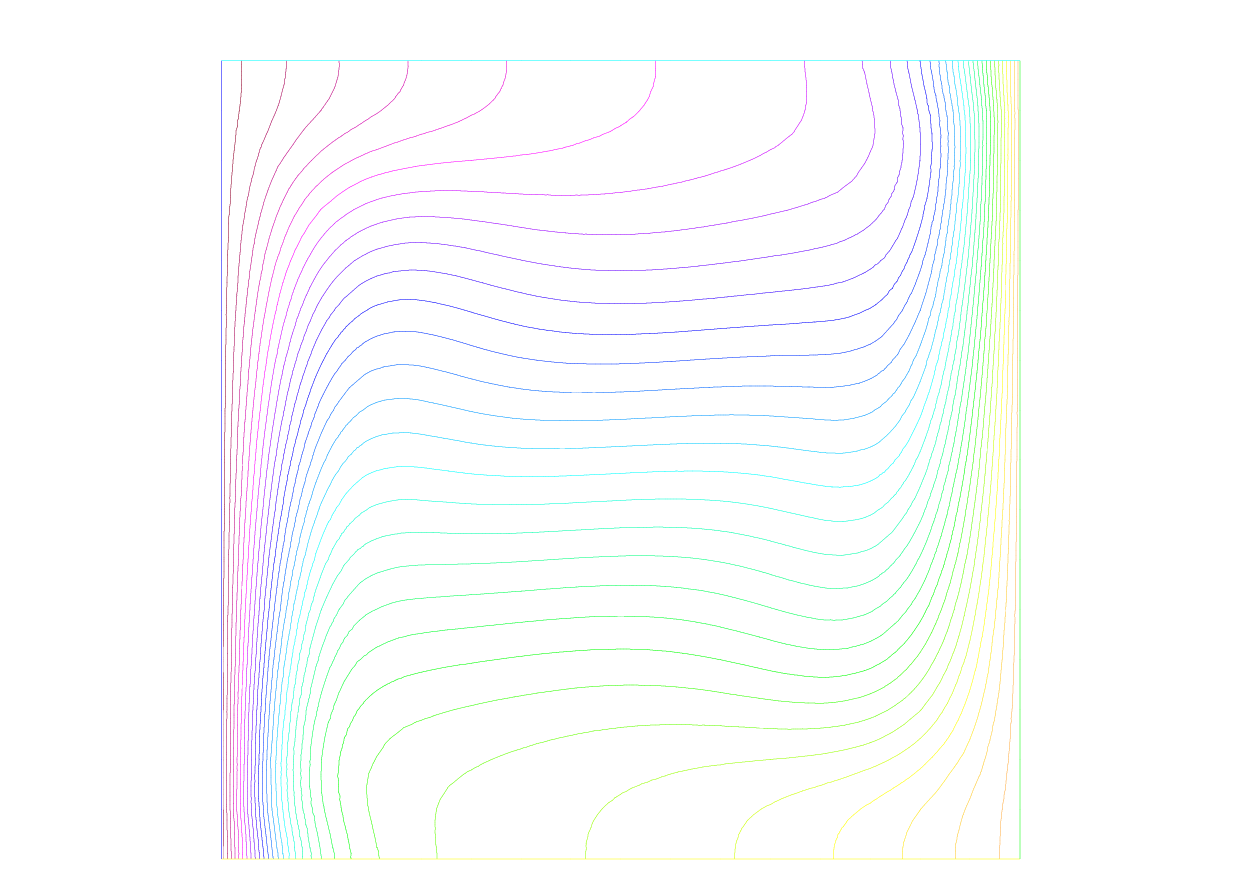}
\caption{Temperature}
\end{subfigure}
\begin{subfigure}{0.3\textwidth}
\includegraphics[width=0.9\linewidth]{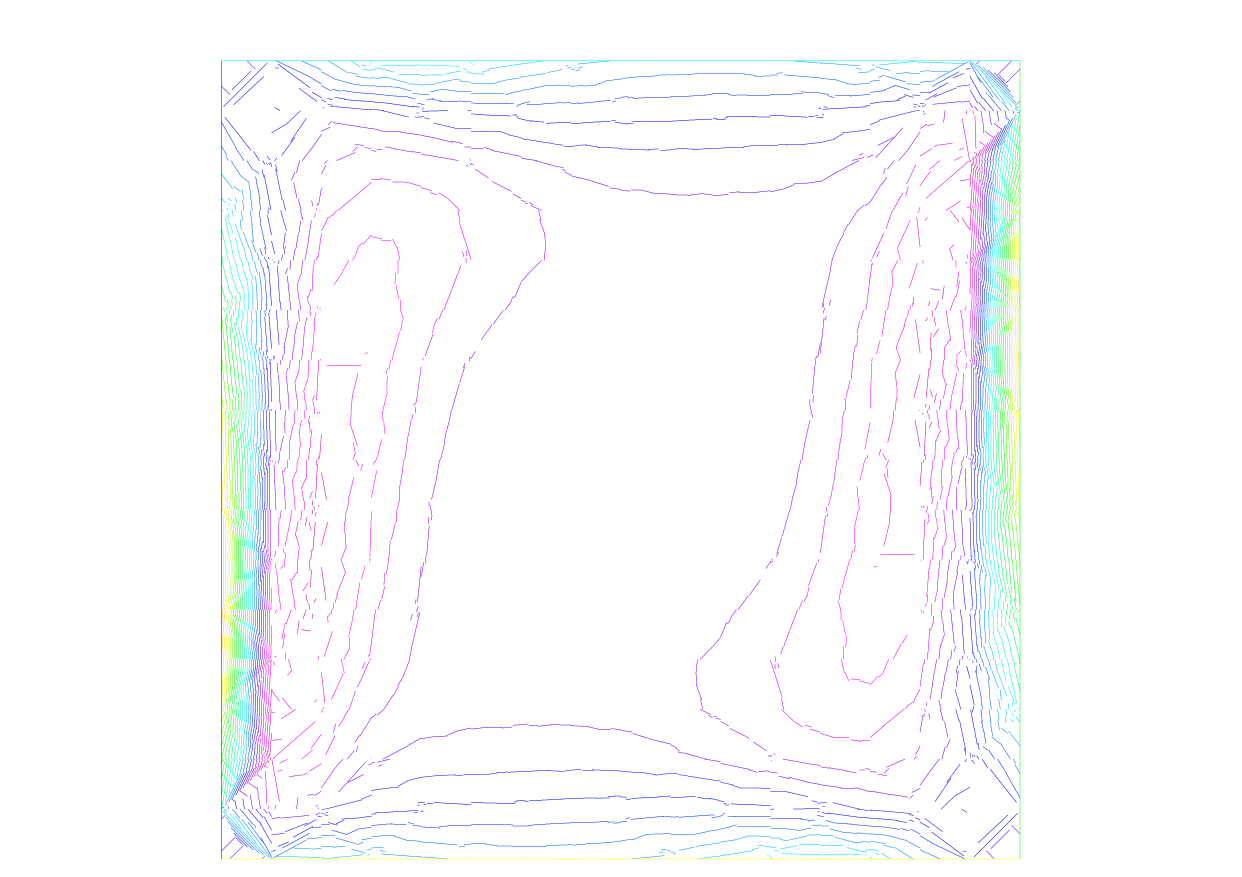}
\caption{Vorticity}
\end{subfigure}
\caption{$Ra=10^5$ Data assimilation solution for $\chi =10^2$.}
\label{fig:image10}
\end{figure}

\begin{figure}[H]
\begin{subfigure}{0.3\textwidth}
\includegraphics[width=0.9\linewidth]{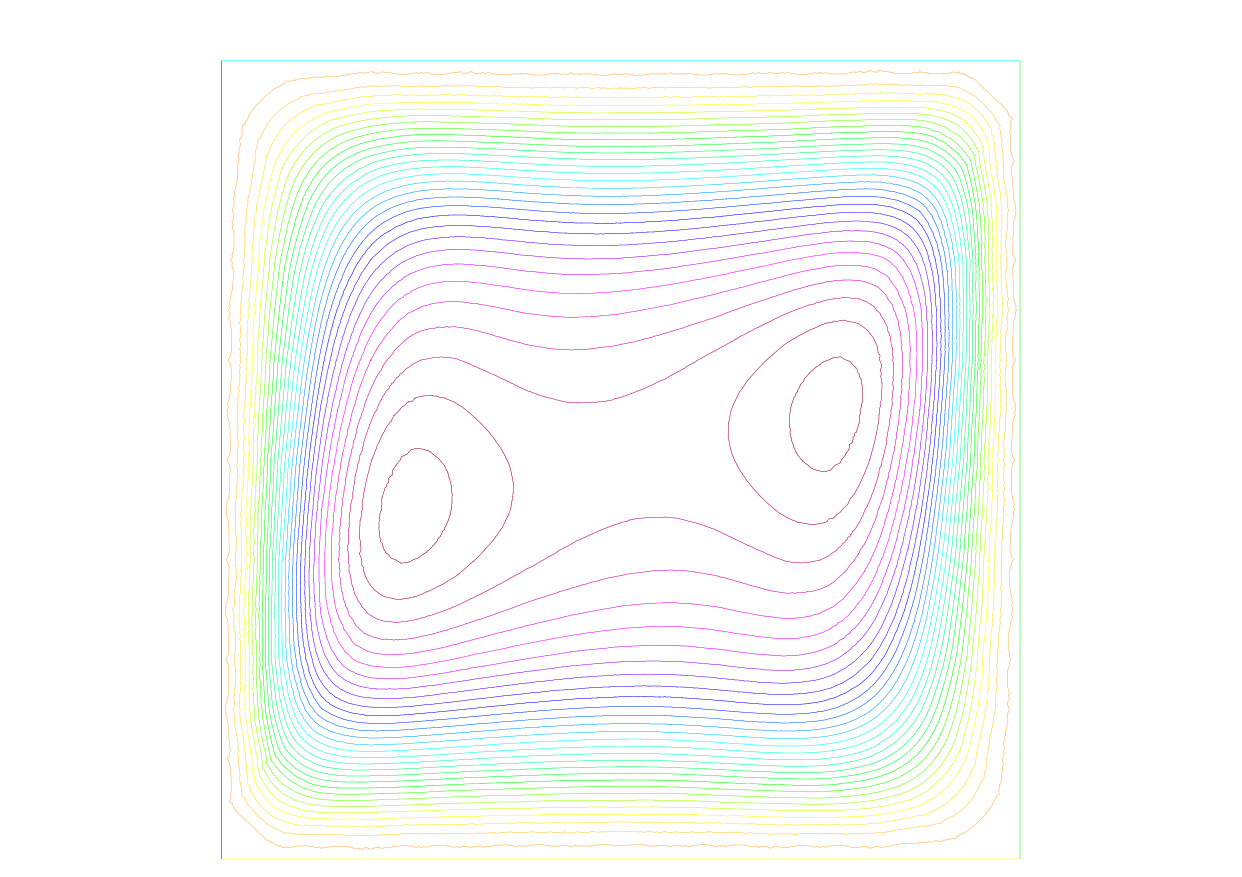} 
\caption{Streamlines}
\end{subfigure}
\begin{subfigure}{0.3\textwidth}
\includegraphics[width=0.9\linewidth]{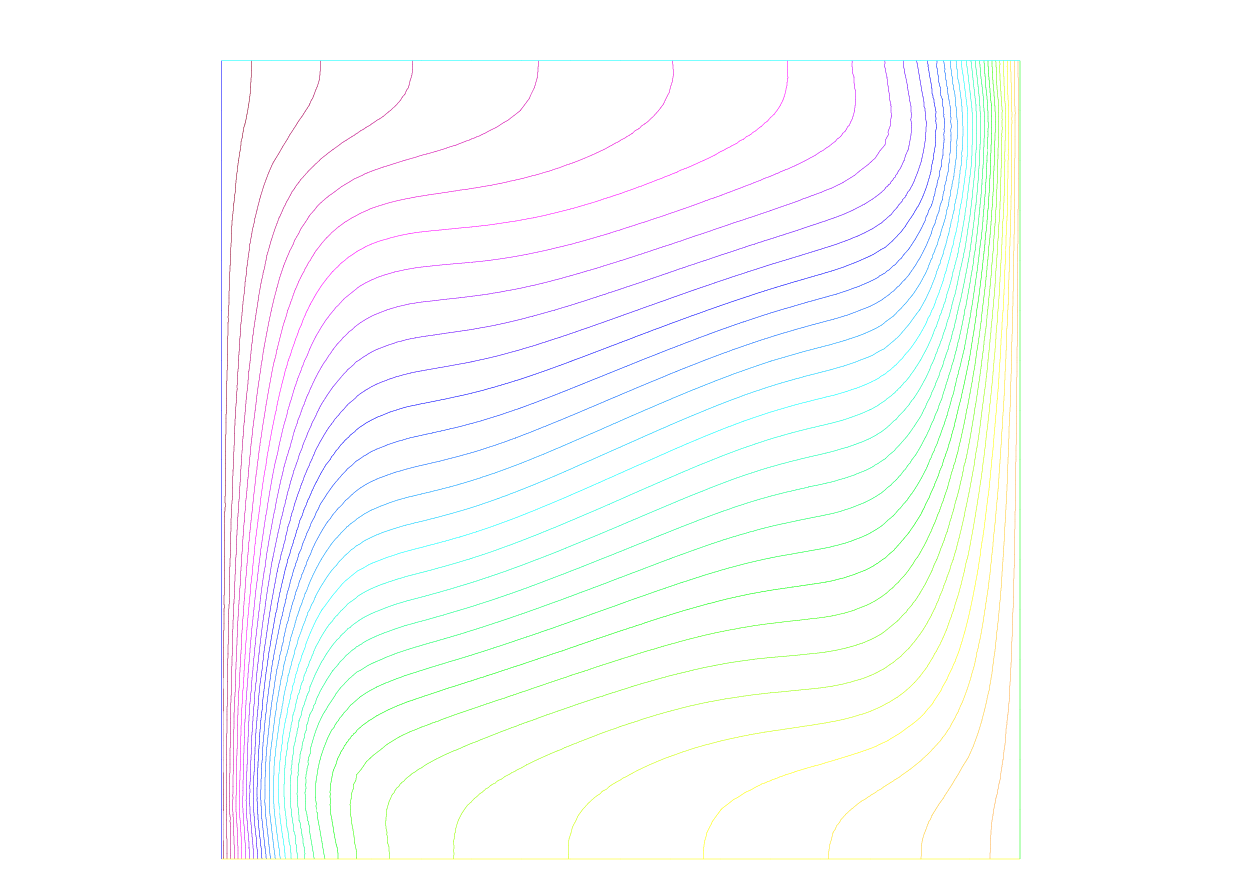}
\caption{Temperature }
\end{subfigure}
\begin{subfigure}{0.3\textwidth}
\includegraphics[width=0.9\linewidth]{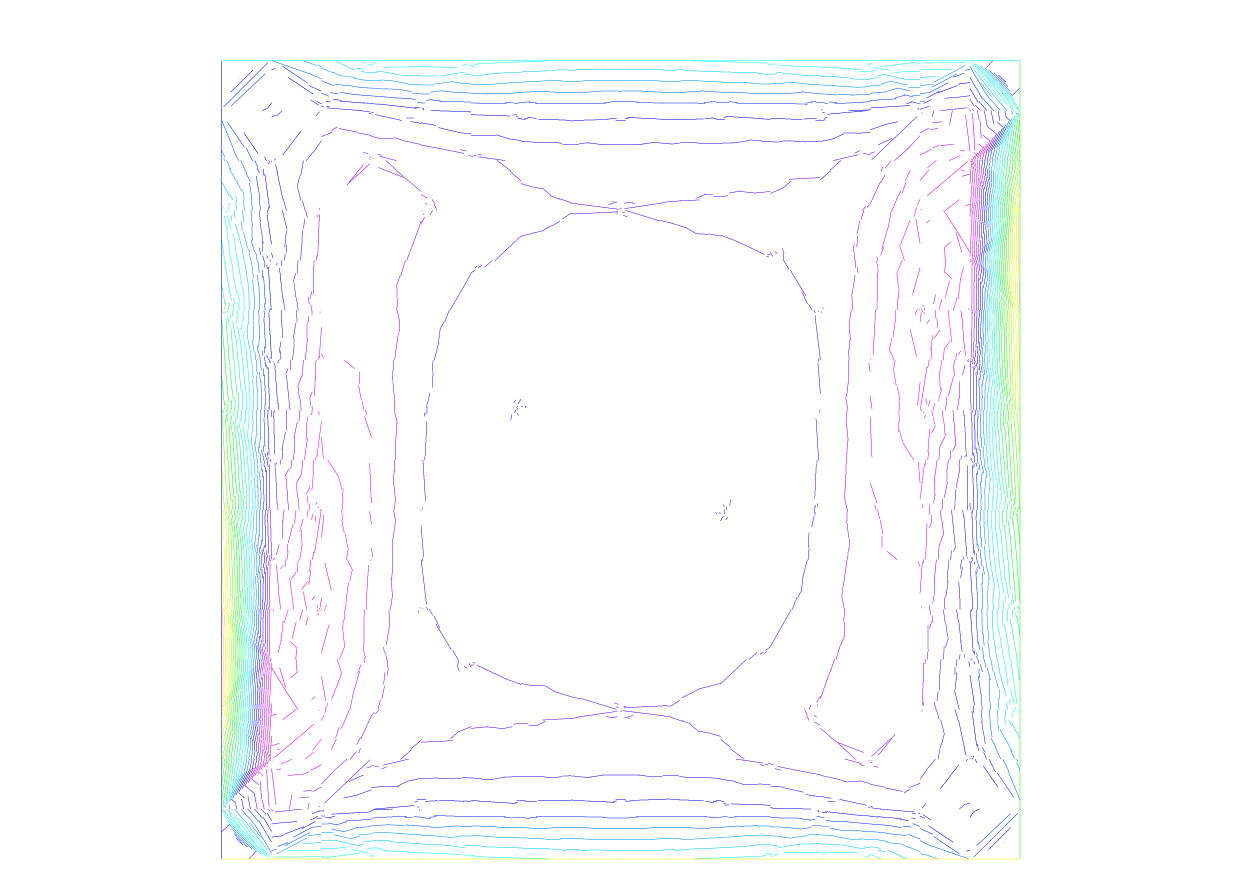}
\caption{Vorticity}
\end{subfigure}
\caption{$Ra=10^5$ Data assimilation solution for $\chi =10^4$.}
\label{fig:image11}
\end{figure}
\begin{figure}[H]
\begin{subfigure}{0.3\textwidth}
\includegraphics[width=0.9\linewidth]{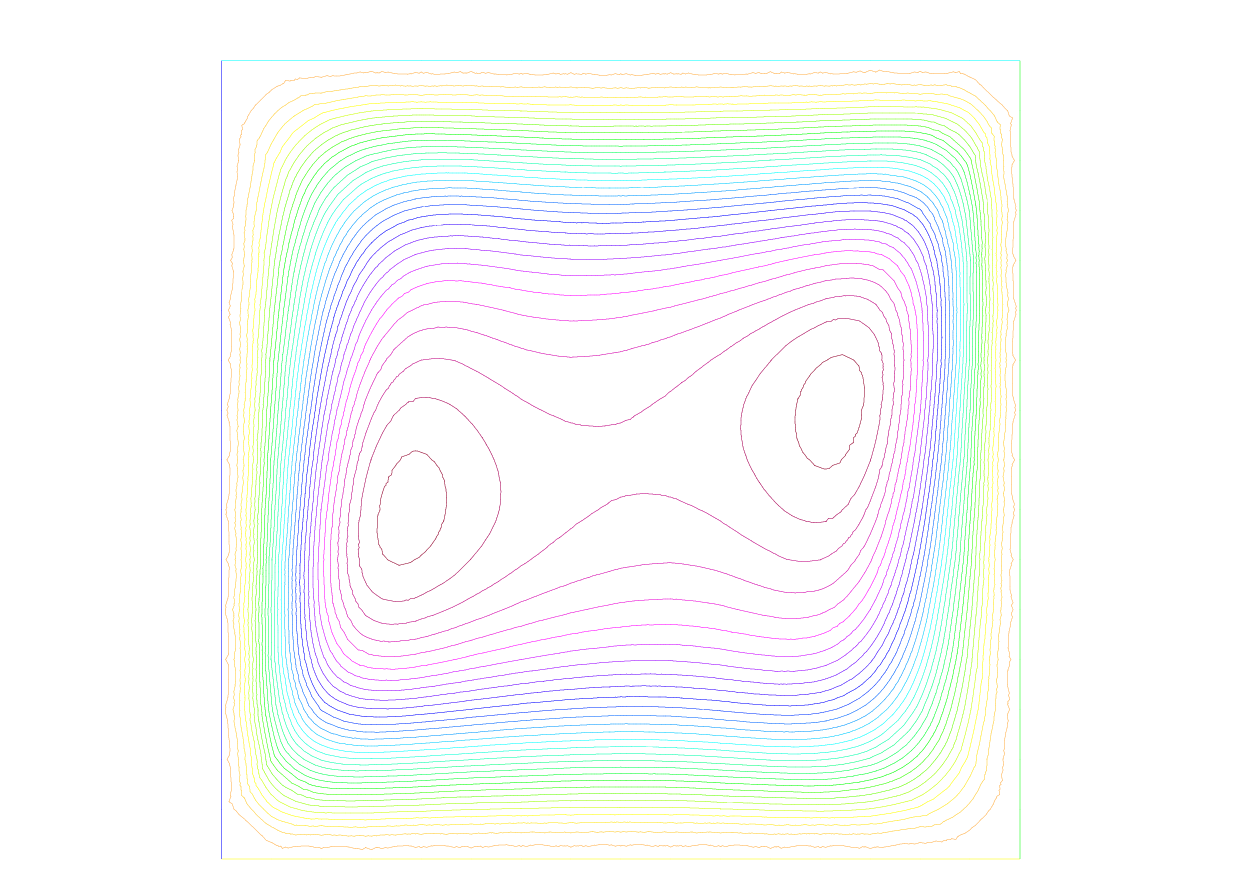} 
\caption{Streamlines}
\end{subfigure}
\begin{subfigure}{0.3\textwidth}
\includegraphics[width=0.9\linewidth]{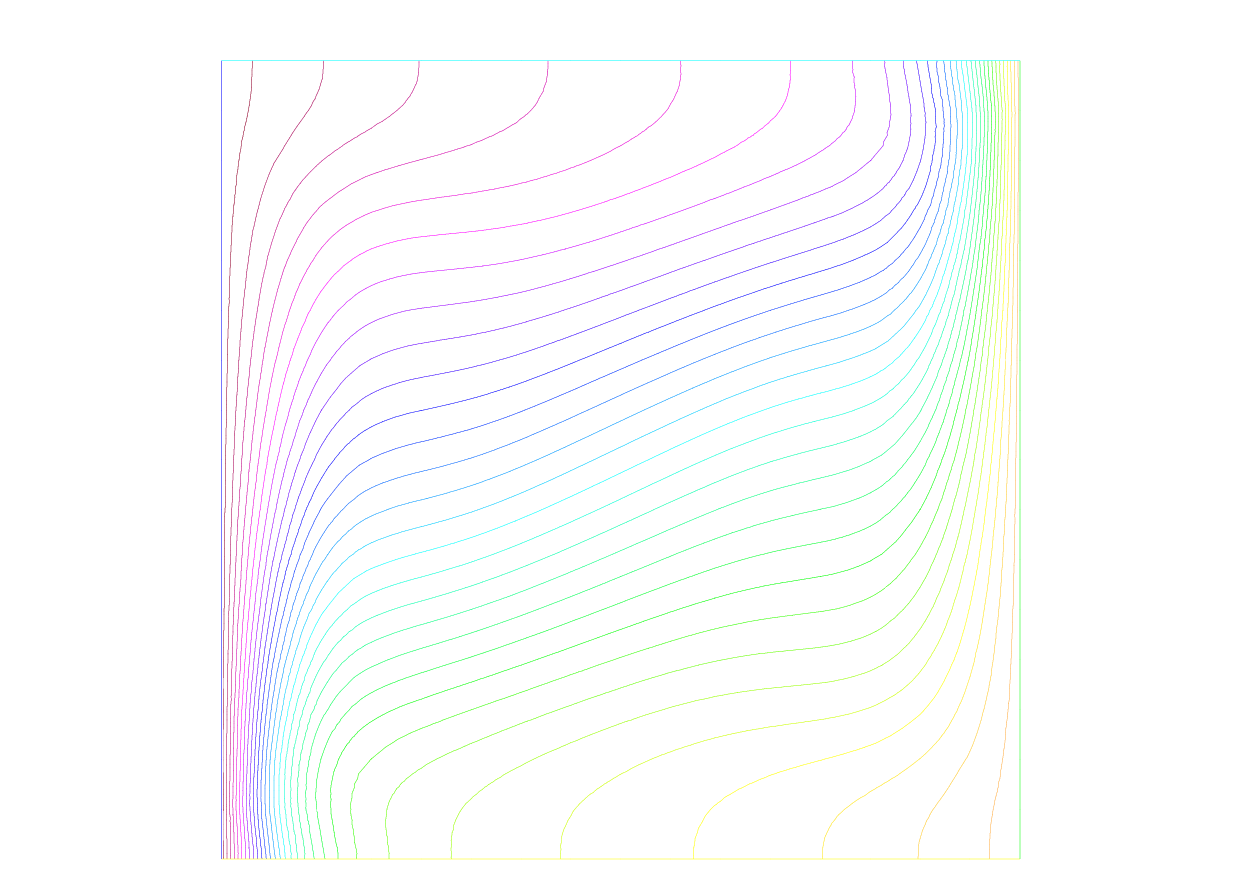}
\caption{Temperature }
\end{subfigure}
\begin{subfigure}{0.3\textwidth}
\includegraphics[width=0.9\linewidth]{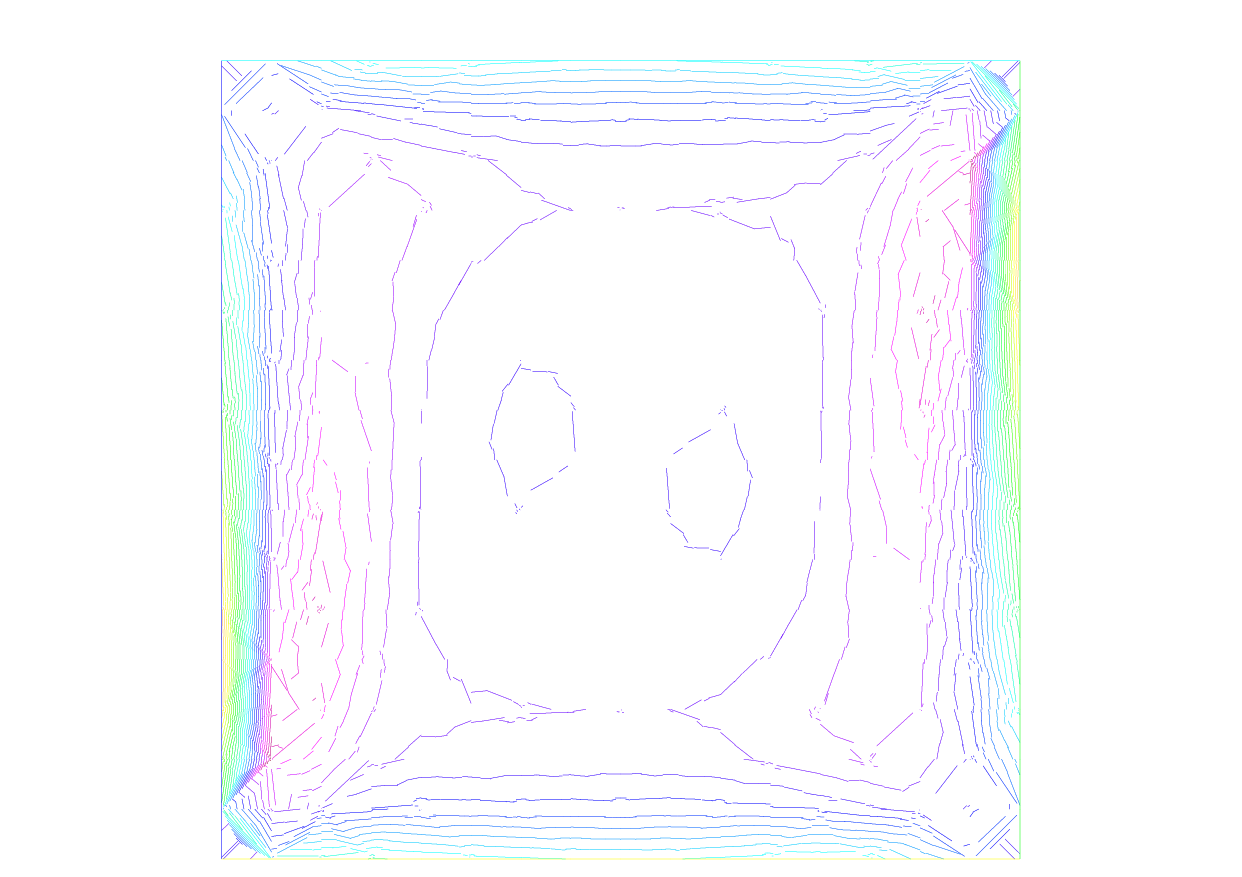}
\caption{Vorticity}
\end{subfigure}
\caption{$Ra=10^5$ Data assimilation solution for $\chi =10^6$.}
\label{fig:image12}
\end{figure}
Figures \ref{fig:image9}-\ref{fig:image12} suggest that increasing $\chi$ matches the solution of (6.1)-(6.3) to the system with Coriolis force. This indicates that the model error has been accounted for correctly.

\begin{figure}[H]
\begin{subfigure}{0.48\textwidth}
\includegraphics[width=0.9\linewidth]{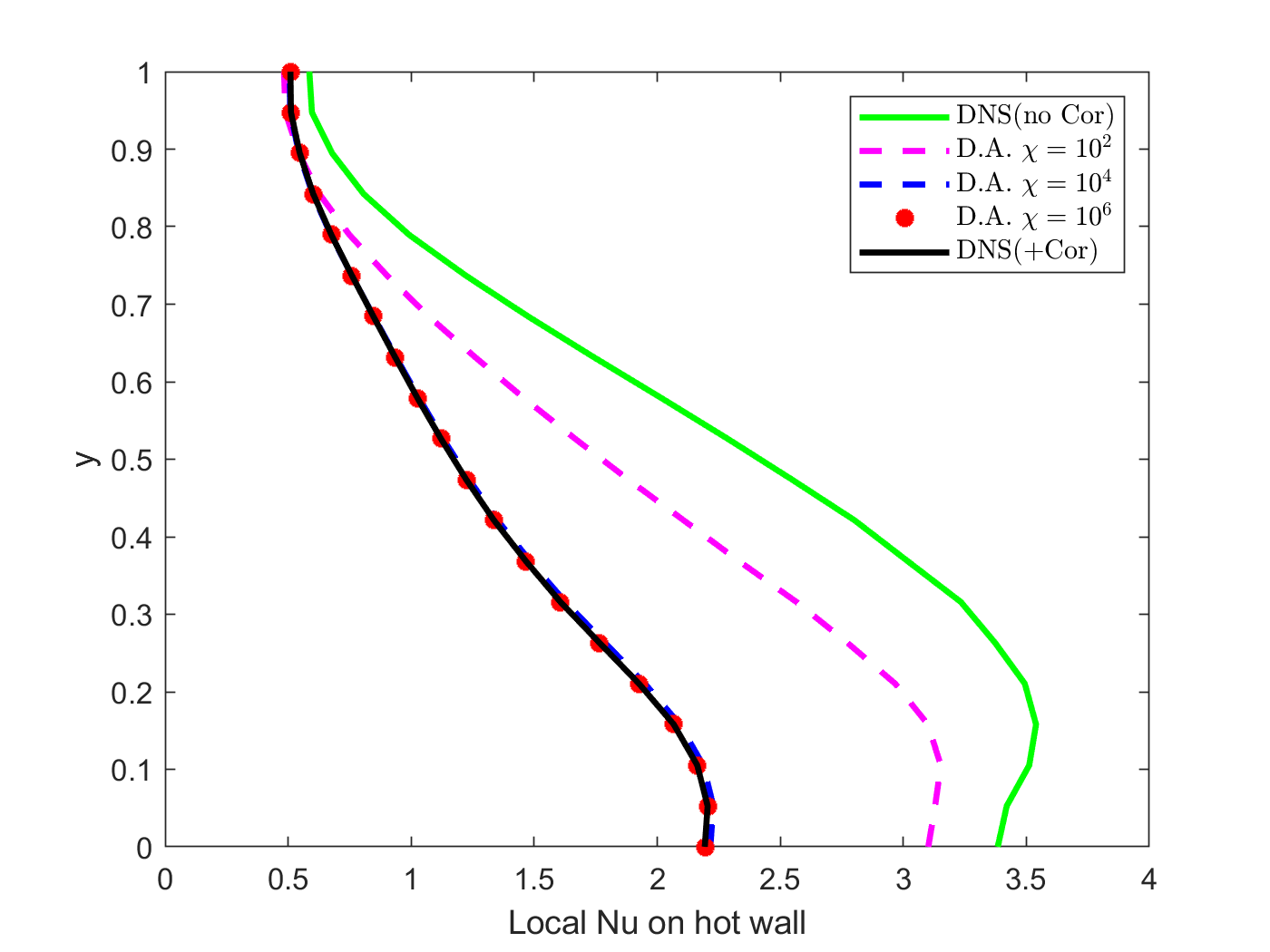} 
\caption{$Ra=10^4$}
\end{subfigure}
\begin{subfigure}{0.48\textwidth}
\includegraphics[width=0.9\linewidth]{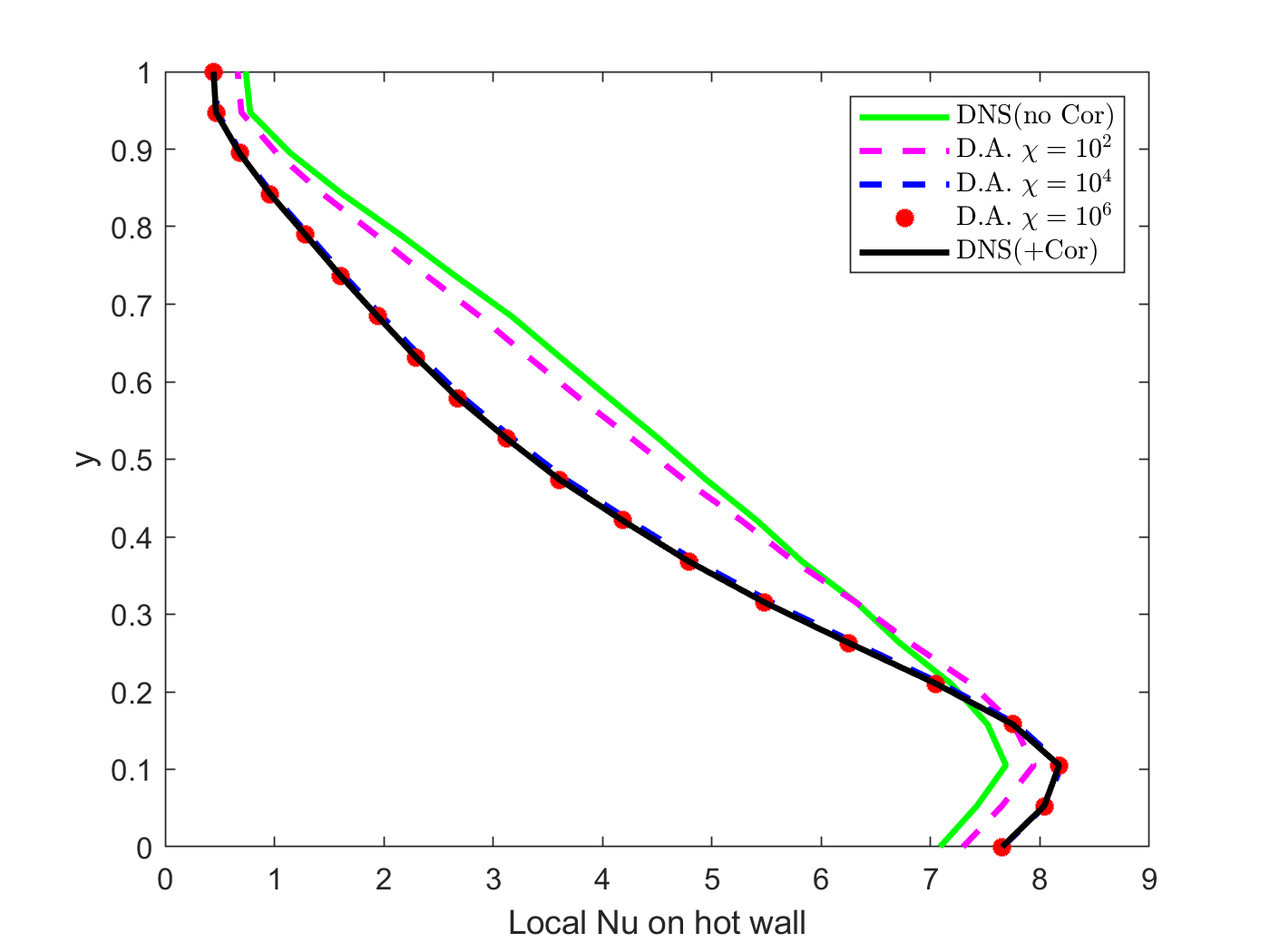}
\caption{$Ra=10^5$ }
\end{subfigure}
\caption{Variation of local Nusselt number on the hot wall for different Rayleigh numbers.}
\label{fig:image113}
\end{figure}

Another quantity of interest for this test problem is the variation of local Nusselt numbers along the vertical walls. $Nu_{local}$ is given by:
$$Nu_{local}=\pm \frac{\partial T}{\partial x},$$
where $\pm$ stands for cold and hot walls, respectively. We calculate and plot $Nu_{local}$ at the hot wall for cases, DNS with and without Coriolis force, and data assimilated solutions for different $\chi$ values, and plot them as given in Figure \ref{fig:image113}.

For both cases, $Ra=10^4$ and $Ra=10^5$, the local Nusselt number at the hot wall approaches the values of DNS solutions with the Coriolis term as $\chi$ increases. These results are also comparable to \cite{cibik11}.
\section{Conclusion}
This report presents an analysis of nudging to correct model errors. Section \ref{stability} establishes the stability of both the true solution 
$u$ and the time-continuous nudged solution $v^h$ in the finite element space. Section \ref{error_estimate_coriolis} provides an error estimate for the nudged solution, comparing it to the true solution $u$, which contains a term omitted from the computational mode $\omega R(u)$. We proved that the error contribution due to model error decays in order of $\mathcal{O}(\chi^{-\frac{1}{2}})$ as $\chi \to \infty$ for both the error without discretization in Theorem \ref {contis_model_error} and continuous in time and discrete in space case. We have verified the expected convergence rate for this method for a common time discretization. Section \ref{temperature_test} considers a natural convection problem with rotation. We have performed comprehensive benchmark tests on a double-pane window problem to study nudging for correcting the model error caused by omitting rotation. Flow patterns of the streamlines, temperature, and vorticity are becoming similar with an increase of $\chi$. We also calculated the local Nusselt number for different $\chi$, and with increased $\chi$, we observed a closer local Nusselt number with the DNS solution with the Coriolis force.  

\section*{Acknowledgment}
This work of William Layton and Rui Fang was partially supported by the NSF grant DMS 2410893. The author Aytekin \c{C}{\i}b{\i}k was partially supported by TUBITAK with the BIDEB-2219 grant.

\bibliography{mybib.bib}
\end{document}